\input amstex
\magnification=1200
\font\cyr=wncyr10
\font\cyb=wncyb10
\font\cyi=wncyi10
\font\cyre=wncyr8
\font\cybe=wncyb8
\documentstyle{amsppt}
\NoRunningHeads
\NoBlackBoxes
\define\sltwo{\operatorname{\frak s\frak l}(2,\Bbb C)}
\define\vir{\operatorname{\frak v\frak i\frak r}}
\define\bnd{\operatorname{\Cal B}}
\define\K{\operatorname{\Cal K}}
\define\HS{\operatorname{\Cal H\Cal S}}
\define\TC{\operatorname{\Cal T\Cal C}}
\define\FR{\operatorname{\Cal F\Cal R}}
\define\md{\operatorname{mod}\ }
\define\mdohbar{\md O(\hbar)}
\define\ad{\operatorname{ad}}
\define\Hilb{\operatorname{Hilb}}
\define\form{\operatorname{form}}
\define\Pontr{\operatorname{Pontr}}
\define\sgn{\operatorname{sgn}}
\define\Vect{\operatorname{Vect}}
\define\jltwo{\operatorname{\frak j\frak l}(2,\Bbb C)}
\define\slstwo{\operatorname{\frak s\frak l}^*(2,\Bbb C)}
\define\End{\operatorname{End}}
\topmatter
\title Approximate representations and Virasoro algebra
\endtitle
\author\bf D.V.Juriev\footnote{\ This is an English translation of the
original Russian version, which is located at the end of the article as an
appendix. In the case of any differences between English and Russian
versions caused by a translation the least has the priority as the original
one.\newline}
\endauthor
\date February, 1997\qquad\qquad\qquad math.RT/9805001\enddate
\affil\eightpoint Research Center for Mathematical Physics and Informatics
``Thalassa Aitheria'',\linebreak
ul.Miklukho-Maklaya 20-180, Moscow 117437 Russia.\linebreak
E-mail: denis\@juriev.msk.ru
\endaffil
\endtopmatter
\document
This article being addressed, in general, to the specialists in the
theory of representations of finite dimensional reductive and infinite
dimensional $\Bbb Z$--graded Lie algebras, the theory of approximations in
functional analysis and asymptotical methods in the operator calculus and,
partially, in the mathematical physics (al\-geb\-ra\-ic quantum field theory),
is devoted to an investigation of algebraic and analytic structure of
infinite dimensional hidden symmetries in the theory of representation of
finite dimensional reductive Lie algebras. More precisely, the subject of
the pa\-per is a realization of the infinite dimensional $\Bbb Z$--graded
Witt algebra of all Lau\-rent polynomial vector fields on a circle and its
one-dimensional nontrivial cent\-ral extension (the Virasoro algebra)
by hidden symmetries in the Verma modules over the Lie algebra $\sltwo$
(the so-called $q_R$--conformal symmetries). These infinite
di\-men\-sional Lie algebras are realized by the $q_R$--conformal symmetries
approximately. In the article there are considered two kinds of
approximate representations: the representations up to a class $\frak S$
of operators (Hilbert-Schmidt operators, compact operators or finite-rank
operators) and asymptotic representations ``$\mdohbar$'', where $\hbar$ is
a parameter (in this case several definitions of an operator ``$O$--large''
by the parameter $\hbar$ are possible). The approximate representations of
the first kind may be naturally considered in the context of the
pseudo-differential calculus [1,2], whereas ones of the second kind -- in
the context of the asymptotic methods [3-5]. Note that the asymptotic
representations ``$\mdohbar$'' were explored in the formalism of the
Karasev-Maslov asymptotic quantisation [6], and representations up to
$\frak S$ are certain generalisations of the projective representations
[7,8].

Results of the article, some of which were anonced in the electronic
preprints of the electronic archive on functional analysis (the Los Alamos
National Laboratories, USA), should be considered more as observations than
as theorems with laborious proofs, which are not clear from their
formulations. Proofs of all facts are performed really in one step and
are bulky computative exercises (omitted where it is possible) based on
the formulas derived by the author earlier; thus, the attention should be
devoted, in general, to the formulations of results themselves.

\head\bf\S1. Preliminary definitions\endhead

\subhead 1.1. The Lie algebra $\sltwo$ and Verma modules over it\endsubhead
Lie algebra $\sltwo$ is a three-dimensional space of $2\times2$ complex
matrices with zero trace supplied by the standard commutator $[X,Y]=XY-YX$,
where the right hand side multiplication is the standard matrix multiplication.
In the basis
$$ l_{-1}=\left(\matrix0&1\\0&0\endmatrix\right),\quad
l_0=\left(\matrix\tfrac12&0\\0&-\tfrac12\endmatrix\right),\quad
l_1=\left(\matrix0&0\\-1&0\endmatrix\right) $$
the commutation relations are of the form: $[l_i,l_j]=(i-j)l_{i+j}$
($i,j=-1,0,1$). Lie algebra $\sltwo$ is $\Bbb Z$--graded:
$\deg(l_i)=-\ad(l_0)l_i=i$, where $\ad(X)$ is the adjoint action operator
in the Lie algebra: $\ad(X)Y=[X,Y]$. Therefore, $\Bbb Z$--graded modules
over $\sltwo$ are $l_0$--diagonal. A vector $v$ in a $\Bbb Z$--graded
module over the Lie algebra $\sltwo$ is called extremal iff $l_1v=0$
and the linear span of vectors $l_{-1}^nv$ ($n\in\Bbb Z_+$) coincides
with the module itself (i.e. $v$ is a cyclic vector). A $\Bbb Z$--graded
module with an extremal vector (in this case it is defined up to a
multiplier) is called extremal [11]. An extremal module is called the
Verma module [12] iff the action of $l_{-1}$ is free in it, i.e. the vectors
$l_{-1}^nv$ are linearly independent. In the case of the Lie algebra $\sltwo$
the Verma modules are just the same as infinite dimensional extremal
modules. An extremal weight of the Verma module is the number defined
by the equality $l_0v=hv$, where $v$ is the extremal vector. The Verma
modules are defined for all complex numbers $h$ and are pairwise nonisomorphic.
Below we shall consider the Verma modules with real extremal weights only.

The Verma module $V_h$ over the Lie algebra $\sltwo$ with the extremal weight
$h$ may be realized in the space $\Bbb C[z]$ of polynomials of a complex
variable $z$. The formulas for the generators of the Lie algebra $\sltwo$ are
of the form:
$$l_{-1}=z,\quad l_0=z\partial_z+h,\quad l_1=z\partial_z^2+2h\partial_z,$$
here $\partial_z=\frac{d}{dz}$.

The Verma module is nondegenerate (i.e. does not contain any proper
sub\-mo\-dule) iff $h\ne-\tfrac{n}2$ ($n\in\Bbb Z_+$). The Verma module $V_h$
is called unitarizable (or hermitean) iff it admits a structure of the
pre-Hilbert space such that $l_i^*=l_{-i}$. The completion of the
unitarizable Verma module will be denoted by $V^{\Hilb}_h$. The Lie algebra
$\sltwo$ acts in $V^{\Hilb}_h$ by the unbounded operators. Also it is rather
useful to consider the formal Verma modules $V^{\form}_h$, which are realized
in the space $\Bbb C[[z]]$ of formal power series of a complex variable $z$,
whereas the formulas for generators of the Lie algebra $\sltwo$ coincide
with ones above. Note that $V_h\subseteq V^{\Hilb}_h\subseteq V^{\form}_h$
and modules $V_h$, $V^{\Hilb}_h$, $V^{\form}_h$ form the Gelfand triple or
the Dirac equipment of the Hilbert space $V^{\Hilb}_h$. An action of the real
form of the Lie algebra $\sltwo$ generated by the anti-Hermitean operators
$il_0$, $l_1-l_{-1}$ and $i(l_1+l_{-1})$ in the
Hilbert space $V^{\Hilb}_h$ by the unbounded operators is exponentiated to
a unitary representation of the corresponding simply connected Lie group.

In the nonunitarizable Verma module over the Lie algebra $\sltwo$ there
exists the unique (up to a scalar multiple) indefinite sesquilinear
form $(\cdot,\cdot)$ such that $(l_iv_1,v_2)=(v_1,l_{-i}v_2)$ for any two
vectors $v_1$ and $v_2$ from the Verma module. If this sesquilinear form
is nondegenerate (in this case the Verma module is nondegenerate) then it
has a signature $(n,\infty)$, where $n$ is finite, and therefore, there is
defined a Pontryagin completion [10] of the Verma module. The corresponding
module in which the Lie algebra $\sltwo$ acts by the unbounded operators
will be denoted by $V^{\Pontr}_h$. The following chain of inclusions holds:
$V_h\subseteq V^{\Pontr}_h\subseteq V^{\form}_h$. An action of the real
form of the Lie algebra $\sltwo$ generated by the anti-Hermitean (with
respect to the nondegenerate indefinite sesquilinear form $(\cdot,\cdot)$)
operators $il_0$, $l_1-l_{-1}$ and $i(l_1+l_{-1})$ in the
Pontryagin space $V^{\Pontr}_h$ by the unbounded operators is exponentiated
to a pseudounitary representation of the corresponding simply connected
Lie group.

\subhead 1.2. Hidden symmetries in Verma modules over the lie algebra
$\sltwo$: Lobachevskii-Berezin $C^*$--algebra and $q_R$--conformal
symmetries\endsubhead

\proclaim{Proposition 1 [13]}\it In the nondegenerate Verma module $V_h$
over the Lie algebra $\sltwo$ there are uniquely defined the operators
$D$ and $F$ such that
$$\aligned [D,l_{-1}]=1,\quad [D,l_0]=D,\quad [D,l_1]=D^2,\\
[l_{-1},F]=1,\quad [l_0,F]=F,\quad [l_1,F]=F^2.\endaligned$$
If the Verma modules are realized in the space $\Bbb C[z]$ of polynomials
of a complex variable $z$ then
$$D=\partial_z,\qquad F=z\tfrac1{\xi+2h},$$
where $\xi=z\partial_z$. The operators $F$ and $D$ obey the following
relations:
$$[FD,DF]=0,\qquad [D,F]=q_R(1-DF)(1-FD),$$
where $q_R=\tfrac1{2h-1}$. In the unitarizable Verma module ($q_R\!\ne\!0$)
the operators $F$ and $D$ are bounded and $F^*=D$, $D^*=F$.
\endproclaim

\rm
The algebra generated by the variables $t$ and $t^*$ with the relations
$[tt^*,t^*t]=0$ and $[t,t^*]=q_R(1-tt^*)(1-t^*t)$ being the Berezin
quantization of the Lobachevskii plane realized in the unit complex disc
(the Poincar\`e realization) [14] is called the Lobachevskii-Berezin algebra.
Proposition 1 allows to consider the Lobachevskii-Berezin algebra as a
$C^*$--algebra. The Lobachevskii-Berezin $C^*$--algebra war re\-dis\-co\-vered
recently by S.Klimek and A.Lesnievsky [15].

\proclaim{Proposition 2 [13]}\it In the nongenerate Verma module $V_h$
over the Lie algebra $\sltwo$ there are uniquely defined the operators $L_n$
($n\in\Bbb Z$) such that
$$[l_i,L_n]=(i-n)L_{i+n}\qquad (i=1,2,1;\ n\in\Bbb Z),$$
moreover, $L_i=l_i$ ($i=-1,0,1$). If the Verma modules are realized in the
space $\Bbb C[z]$ of polynomials of a complex variable $z$ then
$$L_k=(xi+(k+1)h)\partial_z^k\ (k\ge0),\quad
L_{-k}=z^k\tfrac{\xi+(k+1)h}{(\xi+2h)\ldots(\xi+2h+k-1)}\ (k\ge1),$$
where $\xi=z\partial_z$. The operators $L_n$ obey the following relations:
$$[L_n,L_m]=(n-m)L_{n+m},\text{\rm\ if}\ n,m\ge-1\text{\rm\ or}\
n,m\le1.$$
In the unitarizable Verma module the operators $L_n$ are unbounded and
$L_i^*=L_{-1}$.
\endproclaim

\rm
The operators $L_n$ are called the $q_R$--conformal symmetries. They may
be {\it sym\-bo\-li\-cal\-ly\/} represented in the form:
$$L_n=D^{nh}L_0D^{n(1-h)},\qquad L_{-n}=F^{n(1-h)}L_0F^{nh}.$$
To supply the symbolical recording by a sense one should use the general
com\-mu\-tation relations
$$[L_n,f(D)]=(-D)^{n+1}f'(D)\ (n\ge-1),\qquad
[L_{-n},f(F)]=F^{n+1}f'(D)\ (n\ge-1)$$
for $n=0$.

The commutation relations for the operators $D$, $F$ and the generators of
$q_R$--conformal symmetries with the generators of the Lie algebra $\sltwo$
mean that the families $J_k$ and $L_k$ ($k\in\Bbb Z$), where $J_i=D^i$,
$J_{-i}=F^i$ ($i\in\Bbb Z_+$), are families of tensor operators [8,16] for
the Lie algebra $\sltwo$.

To the families of tensor operators $J_k$ and $L_k$ one may correspond
the generating functions -- the operator fields, i.e. the formal Laurent
series of a complex variable $u$ with operator coefficients:
$$J(u)=\sum_{i\in\Bbb Z}J_i(-u)^{-1-i},\quad
T(u)=\sum_{i\in\Bbb Z}L_i(-u)^{-2-i}.$$
For any value of $u$ an operator field realizes a mapping from $V_h$ into
$V^{\form}_h$. The fact that $J_k$ and $T_k$ ($k\in\Bbb Z$) form families
of tensor operators on the language of operator fields means that
$$[l_i,J(u)]=(-u)^iJ(u)-(-u)^{i+1}J'(u),\quad
[l_i,T(u)]=2(-u)^iT(u)-(-u)^{i+1}T'(u).$$
The operator fields $V(u)$, which transform as $s$--differentials under the
action of the Lie algebra $\sltwo$ (i.e. which obey the identity
$[l_i,V(u)]=s(-u)^iV(u)-(-u)^{i+1}V'(u)$), are called the $\sltwo$--primary
fields of spin $s$; so the operator fields $J(u)$ and $T(u)$ are
$\sltwo$--primary ones of spin 1 and 2. The operator fields $J(u)$ and
$T(u)$ as well as their properties are explored in detail in [13].

\subhead 1.3. Infinite dimensional $\Bbb Z$--graded Lie algebras:
the Witt algebra $\frak w^{\Bbb C}$ of all Laurent polynomial vector
fields on a circle and the Virasoro algebra $\vir^{\Bbb C}$, its
one-dimensional nontrivial central extension\endsubhead
The Lie algebra $\Vect(S^1)$ is realized in the space of $C^\infty$--smooth
vector fields $v(t)\partial_t$ on a circle $S^1\simeq\Bbb R/2\pi\Bbb Z$
with the commutator
$$[v_1(t)\partial_t,v_2(t)\partial_t]=
(v_1(t)v_2'(t)-v_1'(t)v_2(t))\partial_t.$$
In the basis $s_n=\sin(nt)\partial_t$, $c_n=\cos(nt)\partial_t$, $h=\partial_t$
the commutation relations have the form:
$$\aligned
[s_n,s_m]&=\tfrac12((m-n)s_{n+m}+\sgn(n-m)(n+m)s_{|n-m|}),\\
[c_n,c_m]&=\tfrac12((n-m)s_{n+m}+\sgn(n-m)(n+m)s_{|n-m|}),\\
[s_n,c_m]&=\tfrac12((m-n)c_{n+m}-(m+n)c_{|n-m|})-n\delta_{nm}h,\\
[h,s_n]&=nc_n,\quad [h,c_n]=-ns_n.
\endaligned$$
Let us denote by $\Vect^{\Bbb C}(S^1)$ the complexification of the Lie algebra
$\Vect(S^1)$. In the basis $e_n=ie^{ikt}\partial_t$ the commutation relations
in the Lie algebra Li $\Vect^{\Bbb C}(S^1)$ have the form:
$$[e_j,e_k]=(j-k)e_{j+k}.$$

It is rather convenient to consider an imbedding of the circle $S^1$ into
the complex plane $\Bbb C$ with the coordinate $z$, so that $z=e^{it}$ on the
circle and the elements of the basis $e_k$ ($k\in\Bbb Z$) are represented by
the Laurent polynomial vector fields: $e_k=z^{k+1}\partial_z$. The
$\Bbb Z$--graded Lie algebra generated by the Laurent polynomial vector
fields (i.e. by the finite linear combinations of elements of the basis
$e_k$) is called the Witt algebra and is denoted by $\frak w^{\Bbb C}$.
The Witt algebra $\frak w^{\Bbb C}$ is the complexification of the
subalgebra $\frak w$ of the algebra $\Vect(S^1)$ generated by the
trigonometric polynomial vector fields on a circle $S^1$, i.e. by the finite
linear combinations of elements of the basis $s_n$, $c_n$ and $h$.

The Lie algebra $\Vect(S^1)$ admits a nontrivial one-dimensional central
extension defined by the Gelfand-Fuchs 2-cocycle [17]:
$$c(v_1(t)\partial_t,v_2(t)\partial_t)=
\int_0^{2\pi}(v_1'(t)v_2''(t)-v_2'(t)v_1'(t))\,dt.$$
This extension being continued to the complexification $\Vect^{\Bbb C}(S^1)$
of the Lie algebra $\Vect(S^1)$ and reduced to the subalgebra
$\frak w^{\Bbb C}$ defines a central one-dimensional ex\-tension of the Witt
algebra, which is called the Virasoro algebra and is denoted by
$\vir^{\Bbb C}$. The Virasoro algebra is generated by the elements
$e_k$ ($k\in\Bbb Z$) and the central element $c$ with the commutation
relations:
$$[e_j,e_k]=(j-k)e_{j+k}+\frac{j^3-j}{12}c$$
and is the complexification of a central extension $\vir$ of the Lie algebra
$\frak w$. In the irreducible representation the central element $c$ of the
Virasoro algebra is mapped to a scalar operator, which is proportional
to the identity operator with a coefficient called {\it the central charge}.

\head\bf\S2. $\HS$--projective representation of the Witt algebra in
Verma modules over the Lie algebra $\sltwo$\endhead

\subhead 2.1. $\frak A$--projective representations [9]\endsubhead

\definition{Definition 1A} Let $\frak A$ be an associative algebra represented
by linear operators in the linear space $H$ and $\frak g$ be the Lie algebra.
The linear mapping $T:\frak g\mapsto\End(H)$ is called {\it a
$\frak A$--projective representation} iff for any $X$ and $Y$ from $\frak g$
there exists an element of the algebra $\frak A$ represented by the
operator $A_{XY}$ such that
$$[T(X),T(Y)]-T([X,Y])=A_{XY}.$$
\enddefinition

If $H$ is infinite dimensional then the representation may be realized
by the unbounded operators with suitable domains of definition.

\remark{Remark 1} Definition 1A may be generalized for arbitrary
anticommutative al\-geb\-ras. In this situation it is deeply related to
the constructions of representations of the anti-commutative algebras
$\jltwo$ and $\slstwo$ from [18,19;20:\S2]. In general, it should be
considered in the context of old ideas of A.I.Maltsev on representations
of arbitarary nonassociative algebras [21]. General anti-commutative
algebras and constructions of their $\frak A$--projective representations
are interesting for the theory of quasi-Hopf algebras, which are
nonassociative as coalgebras, jacobian and c-jacobian quasi-bialgebras
and related structures (see refs in [9]).
\endremark

\remark{Example} If $(\frak g,\frak h)$ is a reductive pair then
any representation of the Lie algebra $\frak g$ is a $\Cal
U(\frak h)$--projective representation of the binary anticommutative algebra
$\frak p$ ($\frak g=\frak h\oplus\frak p$, the binary operation in $\frak p$
has the standard form: $[X,Y]_{\frak p}=\pi([X,Y])$, where $[\cdot,\cdot]$
is the commutator in $\frak g$ and $\pi$ is the projector of $\frak g$
onto $\frak p$ along $\frak h$ [22]).
\endremark

\remark{Remark 2} The standard projective representation is a particular
case of Def.1A if the algebra $\frak A$ is one-dimensional and acts in
$H$ by scalar operators.
\endremark

\remark{Remark 3} If $H$ is a Hilbert (or pre-Hilbert) space then
one may consider the algebra $\HS$ of all Hilbert-Schmidt operators
as $\frak A$. It is possible to consider the algebras $\bnd$, $\K$, $\TC$
and $\FR$ of all bounded, compact, trace class and finite-rank operators,
too.
\endremark

\definition{Definition 1B} Let $\frak A$ be an associative algebra with an
involution $*$ symmetrically represented in the Hilbert space $H$. If
$\frak g$ is a Lie algebra with an involution $*$ then its
$\frak A$--projective representation $T$ in the space $H$ is called
{\it symmetric\/} iff for all elements $a$ from $\frak g$ $T(a^*)=T^*(a)$.
Let $\frak g$ be $\Bbb Z$--graded Lie algebra
($\frak g=\oplus_{n\in\Bbb Z}\frak g_n$) with an involution $*$ such that
$\frak g^*_n=\frak g_{-n}$ and involution is identical on the subalgebra
$\frak g_0$. Let us extend the $\Bbb Z$--grading and the involution $*$
from the Lie algebra $\frak g$ to the tensor algebra $\bold
T^{\cdot}(\frak g)$. A symmetric $\frak A$--projective representation
of $\frak g$ ia called {\it absolutely symmetric\/} iff for any element
$a$ of the algebra $\bold T^{\cdot}(\frak g)$ such that $\deg(a)=0$
the identity $T(a)=T^*(a)$ holds (here the representation $T$ of the
Lie algebra $\frak g$ in $H$ is extended to the mapping of $\bold
T^{\cdot}(\frak g)$ into $\End(H)$.
\enddefinition

\definition{Definition 1C} A $\frak A$--projective representation $T$ of the
Lie algebra $\frak g$ in the linear space $H$ is called {\it almost
absolutely closed\/} iff for any natural $n$ and arbitrary elements
$X_0, X_1, X_2,\ldots X_{n+1}$ of the Lie algebra $\frak g$ there exists
an element $\varphi(X_0,X_1,X_2,\ldots\mathbreak X_{n+1})$ of the algebra
$\frak g$ such that
$$[\ \ldots[[T(X_0),T(X_1)],T(X_2)],\,\ldots,T(X_{n+1})]\!\equiv\!
T(\varphi(X_0,X_1,X_2,\ldots X_{n+1}))\!\!\!\!\pmod{\frak A},$$
here $\frak A$ is considered as represented in $\End(H)$. An almost
absolutely closed $\frak A$--projective representation $T$ of the Lie algebra
$\frak g$ in the linear space $H$ is called {\it absolutely closed\/} iff
$\varphi(\cdot,\ldots,\cdot)\equiv 0$.
\enddefinition

The mappings
$(X_0,X_1,X_2,\ldots X_{n+1})\mapsto\varphi(X_0,X_1,X_2,\ldots X_{n+1})$
associated with an almost absolutely closed $\frak A$--projective
representation of the Lie algebra $\frak g$ define the higher brackets
in the Lie algebra $\frak g$. Objects with higher brackets systematically
appear in many branches of mathematics and mathematical physics (see e.g.
the book [6] or the article [23] among others and numerous references
wherein).

\subhead 2.2. $\HS$--projective representation of the Witt algebra by
$q_R$--conformal symmetries in the unitarizable Verma modules $V_h$ over
the Lie algebra $\sltwo$ [9]
\endsubhead
Note that the Witt algebra $\frak w^{\Bbb C}$ admits a natural involution $*$.

\proclaim{Theorem 1A}\it The generators $L_k$ ($k\in\Bbb Z$) of
$q_R$--conformal symmetries in the unitarizable Verma module $V_h$ over the
Lie algebra $\sltwo$ realize an absolutely symmetric $\HS$--projective
representation of the Witt algebra $\frak w^{\Bbb C}$. Adding the tensor
operators $J_k$ ($k\in\Bbb Z$) to the tensor operators $L_k$ one receives
an absolutely symmetric $\HS$--projective representation of a semi-direct
sum of the Witt algebra and the infinite-dimensional $\Bbb Z$--graded
Heisenberg algebra (the one-dimensional central ex\-tension of the
infinite-dimensional $\Bbb Z$--graded abelian Lie algebra $\Bbb C[z,z^{-1}]$
of Laurent polynomials).
\endproclaim

The statement of the theorem follows from the explicit formulas for
generators of $q_R$--conformal symmetries and tensor operators $J_k$.
A verification on the belonging to the class $\HS$ of
Hilbert-Schmidt operators for operators of the fixed degree with respect
to the grading in the $\Bbb Z$--graded space of polynomials $\Bbb C[z]$
supplied by the certain scalar product with orthogonal weight spaces
(such as any unitarizable Verma module over the Lie algebra $\sltwo$ is)
does not produce any problems.

\remark{Remark 4} The $\HS$--projective representations of the Witt algebra
in the uni\-ta\-ri\-zable Verma modules over $\sltwo$ are absolutely closed.
\endremark

\remark{Remark 5} Theorem 1A is generalized for the pseudo-unitary case
with the sub\-s\-ti\-tution of the class $\HS$ of Hilbert-Schmidt operators
to the class $\K$ of compact operators.
\endremark

The results on the ``exponentiated'' version of Theorem 1A are announced
in the e-prints [24,25].

\subhead 2.3. $\FR$--projective representations of the Witt algebra in
Verma mo\-dules over the Lie algebra $\sltwo$ with the extremal weight
$1$ and $\frac12$\endsubhead
The sattement of Theorem 1A may be strengthened in some particular cases.

\proclaim{Theorem 1B}\it For $h=1$ or $h=\frac12$ the $\HS$--projective
representation of the Witt algebra by $q_R$--conformal symmetries in
the unitarizable Verma modules over the Lie algebra $\sltwo$ is
$\FR$--projective. For $h=\frac12$ the $\HS$--projective representation
of a semi-direct sum of the Witt algebra and the infinite-dimensional
Heisenberg algebra is also $\FR$--projective.
\endproclaim

The statement follows from the explicit formulas for tensor operators
$L_k$ and $J_k$.

\remark{Remark 6} For $h=0$ the $\K$--projective representation of the
Witt algebra in the Verma module over the Lie algebra $\sltwo$ is
$\FR$--projective.
\endremark

\subhead 2.4. $\bnd$--projective representations of the Witt
algebra\endsubhead
Note that one of disadvantages of many $\frak A$--projective
representations, in particular, if $\frak A$ is a certain operator class
$\K$, $\HS$, $\TC$ or $\FR$, is that their set is not closed under tensor
products, in general. For the class $\bnd$ of bounded operators the tensor
products of $\bnd$--projective representations are defined. However, any
$\frak S$--projective representation, where $\frak S$ is one of classes above,
is $\bnd$--projective and as such is sometimes nontrivial (when the
initial representation is realized by unbounded operators), that allows to
construct their tensor products, which will be nontrivial
$\bnd$--projective representations in this case.

Let us call a $\frak S$--projective representation of the Lie algebra
$\frak g$ ($\frak S$ is a certain operator class) {\it
$\frak S$--irreducible\/} iff the representation operators can not
be si\-mul\-ta\-neous\-ly transform by an addition of elements of $\frak S$
to the form, in which they have a common proper invariant subspace.

\proclaim{Theorem 2}\it Let $T_h$ be the $\bnd$--projective representation
of the Witt algebra $\frak w^{\Bbb C}$ by the $q_R$--conformal symmetries
in the Verma module $V_h$ over the Lie algebra $\sltwo$ then the
$\bnd$--projective representations $S^n(T_h)$ (here $S^n$ denotes
the functor of the $n$-th symmetric degree, cf.[26]) are $\bnd$--irreducible.
\endproclaim

The statement of the theorem follows from the explicit formulas for
generators of $q_R$--conformal symmetries.

\remark{Remark 7} Decompositions of the tensor products
$T_{h_1}\otimes\ldots\otimes T_{h_n}$ of $\bnd$--projective representations
$T_{h_i}$ on $\bnd$--irreducible components are not known in general case.
\endremark

\head\bf\S3. Asymptotic representations of the Witt algebra and the
Virasoro algebra in Verma modules over the Lie algebra $\sltwo$\endhead

\subhead 3.1. The asymptotics of $q_R$--conformal symmetries in Verma
modules $V_h$ over the Lie algebra $\sltwo$ at $q_R\to0$ ($h\to\frac12$)
and at $q_R\to1$ ($h\to1$) and the Witt algebra $\frak w^{\Bbb C}$\endsubhead
Let $\frak A$  be the algebra of all finite linear combinations of
expressions of the form $f(qp)p^n$ ($n\in\Bbb Z_+$) or $q^nf(qp)$
($n\in\Bbb Z_+$), where $[p,q]=1$ and $f$ are rational functions, whose
denominators have no zeroes in integer points. Let us call the algebra
$\frak A$ {\it the extended Weyl algebra}. The extended Weyl algebra
is a topological algebra with respect to natural convergence.

\proclaim{Lemma}\it The extended Weyl algebra $\frak A$ admits a strict
representation in the Verma module $V_h$ $(h\not\in\Bbb Z/2)$:
$$p\to\partial_z,\qquad q\to z.$$
Generators of $q_R$--conformal symmetries $L_k$ $(k\in\Bbb Z)$ as well as
generators $D$, $F$ of the Lobachevskii-Berezin algebra belong to the image
of the extended Weyl algebra $\frak A$.
\endproclaim

Below we shall interpret the ``$O$--large'' for the asymptotics of
operators in the Verma modules over the Lie algebra $\sltwo$ with
non-half-integer extremal weights, which belong to the extended Weyl
algebra (because its representations in the Verma modules are strict one
can identify the algebra itself with its image) in the sense of topology in
the extended Weyl algebra.

\proclaim{Theorem 3A}\it Generators of $q_R$--conformal symmetries in the
Verma modules $V_{\frac12+\hbar}$ $(0\!<\!\hbar\!<\!\frac12)$ over the Lie
algebra $\sltwo$ form an asymptotic ``$\mdohbar$'' representation of the
Witt algebra $\frak w^{\Bbb C}$.

Generators of $q_R$--conformal symmetries in the Verma modules $V_{1\pm\hbar}$
$(0\!<\!\hbar\!<\!\frac12)$ over the Lie algebra $\sltwo$ form an asymptotic
``$\mdohbar$'' representation of the Witt algebra $\frak w^{\Bbb C}$, too.
\endproclaim

The first statement of the theorem, which immediately follows from the
explicit form of generators of $q_R$--conformal symmetries, was in fact
proved in [10], where the corresponding ``exponentiated'' version was
considered. The second is similar to the first.

\remark{Remark 8} The statement of the theorem is not true if the
``$O$--large'' is considered in sense of (weak) operator convergence
in the space of operators in the Verma modules $V_h$ over the Lie algebra
$\sltwo$ identified with the space $\Bbb C[z]$ of polynomials.
\endremark

An answer on the question on a ``divergence'' between theorem 3A and the
remark 8 is in the following: the strict representations of the extended
Weyl algebra do not even exist for the limit value $\hbar\!=\!0$, but
being continuous (in a weak operator topology in the space $V_h$) for
all values of the parameter $\hbar$ from an interval $0\!<\!\hbar\!<\!\frac12$,
are not {\it uniformly\/} continuous on this interval. Note that the
continuous representations of the extended Weyl algebra in the verma modules
$V_h$ over the Lie algebra $\sltwo$ can be extended to representations by
unbounded operators in the spaces $V^{\Hilb}_h$ (or $V^{\Pontr}_h$), and,
therefore, uncontinuous, that, however, is typical for the representation
theory of Lie algebras.

\subhead 3.2. The estimation ``$\md\HS$'' of the second term of the
asymptotics of $q_R$--conformal symmetries in Verma modules $V_h$ over
the Lie algebra $\sltwo$ at $q_R\!\to\!0$ ($h\!\to\!\frac12$) and at
$q_R\!\to\!1$ ($h\!\to\!1$) and the Virasoro algebra $\vir^{\Bbb C}$\endsubhead
Under the ``hybridization'' of the approximations ``$\md O(\hbar^n)$'' and
``$\md\frak S$'' an interesting phenomenon of ``noncommutativity'' of
estimations is observed, the asymptotics of $q_R$--conformal symmetries
supply us by an example.

\proclaim{Theorem 3B} For generators of $q_R$--conformal symmetries in the
Verma modules $V_{\frac12+\hbar}$ $(0\!<\!\hbar\!<\!\frac12)$ over the Lie
algebra $\sltwo$ the following asymptotics ``$\md O(\hbar^2)$'' exists:
$$[L_i,L_j]=(i-j)L_{i+j}+\hbar\Cal A+O(\hbar^2),$$
where
$$\Cal A\equiv\frac23(i^3-i)\delta_{ij} \mod{\HS}.$$
So generators of $q_R$--conformal symmetries form an approximate representation
of the Virasoro algebra $\vir^{\Bbb C}$ for $h\to\frac12$ in the specified
sense.

For generators of $q_R$--conformal symmetries in the Verma modules
$V_{1\pm\hbar}$ $(0\!<\!\hbar\!<\!\frac12)$ over the Lie algebra $\sltwo$
the following asymptotics ``$\md O(\hbar^2)$'' exists, too:
$$[L_i,L_j]=(i-j)L_{i+j}+\hbar\Cal A+O(\hbar^2),$$
where
$$\Cal A\equiv\frac23(i^3-i)\delta_{ij} \mod{\HS}.$$
So generators of $q_R$--conformal symmetries form an approximate representation
of the Virasoro algebra $\vir^{\Bbb C}$ for $h\to1$ in the specified sense.
\endproclaim

In the statement of the theorem the estimations ``$\md O(\hbar^2)$''
are interpreted in the sense of a convergence in the extended Weyl algebra
$\frak A$, and $\HS$ denotes the set of all elements of this algebra, which
are realized by Hilbert-Schmidt operators under the strict representation
in the Verma modules over the Lie algebra $\sltwo$ for all sufficiently
small values of $\hbar$ (otherwords, if $\pi_{h}$ is the representation of
$\frak A$ in $V_h$ and $\HS_h$ is the space of Hilbert-Schmidt operators in
$V_h$ then in the theorem $\HS$ is interpreted as
$$\lim_{\hbar\to0}\bigcap_{h_0<h<h_0+\hbar}\pi_{h}^{-1}(\HS_{h}\cap
\pi_{h}(\frak A)),$$
where $h_0=\tfrac12$ or $h_0=1$).

\demo\nofrills{} The result of Theorem 3B follows from explicit computations.
Let us compute the commutator $[L_2,L_{-2}]$ and estimate it that is enough
for the determination of the central charge of the Virasoro algebra.
The Verma module $V_h$ is realized in the space $\Bbb C[z]$ of polynomials,
$\xi=z\partial_z$, and generators of $q_R$--conformal symmetries are
defined by expressions written above. Then for $h\to\frac12$
$$\aligned
[L_2,L_{-2}]=&\tfrac{(\xi+3h)^2(\xi+1)(\xi+2)}{(\xi+2h)(\xi+2h+1)}-
\tfrac{(\xi+3h-2)^2\xi(\xi-1)}{(\xi+2h-1)(\xi+2h-2)}=\\
&\tfrac{(\xi+\tfrac32+3\hbar)^2(\xi+1)(\xi+2)}{(\xi+1+2\hbar)(\xi+2+2\hbar)}-
\tfrac{(\xi-\tfrac12+3\hbar)^2\xi(\xi-1)}{(\xi-1+2\hbar)(\xi+2\hbar)}.
\endaligned$$
The following estimation ``$\md O(\hbar^2)$'' holds ($\hbar$ in denominators
is not excluded because for the zero value of $\hbar$ a representation
of the extended Weyl algebra in the corresponding Verma module is not defined):
$$\aligned
[L_2,L_{-2}]\equiv &(\xi+\tfrac32+3\hbar)^2(1-\tfrac{2\hbar}{\xi+1+\hbar})
(1-\tfrac{2\hbar}{\xi+2+\hbar})-\\
&(\xi-\tfrac12+3\hbar)^2(1-\tfrac{2\hbar}{\xi-1+\hbar})
(1-\tfrac{2\hbar}{\xi+\hbar})\mod{O(\hbar^2)}.
\endaligned$$
Extracting the terms of order $\hbar$ (because the first term of the order 1
is known) one receives:
$$2\hbar(\xi-\tfrac12)^2(\tfrac1{\xi-1+\hbar}+\tfrac1{\xi+\hbar})-
2\hbar(\xi+\tfrac32)^2(\tfrac1{\xi+1+\hbar}+\tfrac1{\xi+2+\hbar})+12\hbar.$$
Perform the estimation ``$\md\HS$'' and obtain:
$$\aligned
&2\hbar\xi^2(\tfrac1{\xi-1+\hbar}+\tfrac1{\xi+\hbar}-\tfrac1{\xi+1+\hbar}-
\tfrac1{\xi+2+\hbar})-12\hbar-\\
&2\hbar\xi(\tfrac1{\xi-1-\hbar}+\tfrac1{\xi+\hbar})-
6\hbar(\tfrac1{\xi+1+\hbar}+\tfrac1{\xi+2+\hbar})\sim\\
&2\hbar\xi^2(\tfrac2{(\xi-1+\hbar)(\xi+1+\hbar)}+
\tfrac2{(\xi+\hbar)(\xi+2+\hbar)})+12\hbar-16\hbar\sim\\
&8\hbar+12\hbar-16\hbar=4\hbar.
\endaligned$$
So the ``asymptotic'' central charge is equal to $8\hbar$.

Analogously, for $h\to1$
$$\aligned
[L_2,L_{-2}]=&\tfrac{(\xi+3h)^2(\xi+1)(\xi+2)}{(\xi+2h)(\xi+2h+1)}-
\tfrac{(\xi+3h-2)^2\xi(\xi-1)}{(\xi+2h-1)(\xi+2h-2)}=\\
&\tfrac{(\xi+3+3\hbar)^2(\xi+1)(\xi+2)}{(\xi+2+2\hbar)(\xi+3+2\hbar)}-
\tfrac{(\xi+1+3\hbar)^2\xi(\xi-1)}{(\xi+2\hbar)(\xi+1+2\hbar)}.
\endaligned$$
The following estimation ``$\md O(\hbar^2)$'' holds ($\hbar$ in denominators
is not excluded because for the zero value of $\hbar$ a representation
of the extended Weyl algebra in the corresponding Verma module is not defined):
$$\aligned
[L_2,L_{-2}]\equiv &(\xi+1)(\xi+3+3\hbar)(1-\tfrac{2\hbar}{\xi+2+\hbar})
(1-\tfrac{2\hbar}{\xi+3+\hbar})(1+\tfrac{3\hbar}{\xi+3+\hbar})-\\
&(\xi-1)(\xi+1+3\hbar)(1-\tfrac{2\hbar}{\xi+\hbar})
(1-\tfrac{2\hbar}{\xi+1+\hbar})(1+\tfrac{3\hbar}{\xi+1+\hbar})\mod{O(\hbar^2)}.
\endaligned$$
Extracting the terms of order $\hbar$ (because the first term of the order 1
is known) one receives:
$$\hbar(\xi+1)(\xi+3)(\tfrac1{\xi+3+\hbar}-\tfrac2{\xi+2+\hbar})-
\hbar(\xi+1)(\xi-1)(\tfrac1{\xi+1+\hbar}-\tfrac2{\xi+\hbar})+6\hbar.$$
Perform the estimation ``$\md\HS$'' and obtain:
$$8\hbar+2(\xi+1)(\tfrac{\xi-1}{\xi+\hbar}-\tfrac{\xi+3}{\xi+2+\hbar})\sim
8\hbar-4\hbar=4\hbar.$$
So the ``asymptotic'' central charge is equal to $8\hbar$.
\enddemo

\remark{Remark 9} The estimation ``$\md O(\hbar^2)$'' may be strengthened
to ``$\md O(\hbar^n)$'' for any {\it finite\/} $n$ (but not to
``$\md O(\hbar^\infty)$'' in view of the further estimation ``$\md\HS$'').
\endremark

\remark{Remark 10} As follows from the results of the second paragraph
the change of an order of estimations gives the zero central charge of the
Virasoro algebra.
\endremark

Perhaps, the noncommutativity of an order of the estimations
``$\md O(\hbar^n)$'' and ``$\md\HS$'' in the case of
$q_R$--conformal symmetries is a partial case of a more general and
fundamental fact of deviations between asymptotical theory of
pseudodifferential operators and pseudodifferential calculus on
asymptotic manifolds in sense of [27].

Note that the result of Theorem 3B contains a strange coincidence
of ``asymptotic'' central charges of the Virasoro algebra for
$h\to\frac12$ and $h\to1$ (namely, $c=8\hbar$, $\hbar\to0$).
Perhaps, this coincidence may be explaned by the hypothetical suggestion
that for all values of the extremal weight $h$ the Virasoro algebra
is presented on a certain more hidden and unexplicated level universally
(i.e. its characteristics do not depend on $h$).

\head\bf Conclusion
\endhead

In the paper there are investigated various approximate representations of
the infinite dimensional $\Bbb Z$--graded Lie algebras: the Witt algebra
of all Laurent po\-ly\-no\-mial vector fields on a circle and its
one-dimensional nontrivial central extension, the Virasoro algebra, by the
infinite dimensional hidden symmetries in the Verma modules over the Lie
algebra $\sltwo$. There are considered as asymptotic rep\-re\-sen\-tations
``$\md O(\hbar^n)$'' and representations up to a class $\frak S$ of operators
(compact operators, Hilbert-Schmidt operators and finite-rank operators) as
cases, which combine both types of approximations (and in these cases an
effect of non\-com\-mu\-ta\-ti\-vi\-ty of the order of their realizations
is explicated, that perhaps is underlied by a more general fundamental fact
of deviations between the asymptotical theory of pseudodifferential operators
and the pseudodifferential calculus on the asymptotic manifolds in sense of
[27]). Some applications of the discussed problems to the applied problems
of information technologies (an organization of the information transmission
in integrated videocognitive interactive systems for accelerated non\-ver\-bal
computer and telecommunications) were investigated in the report [28].

\ \newline
\ \newline
\ \newline

\centerline{\bf References}
\eightpoint
\roster
\item" [1]" Taylor M., Pseudo-differential operators. B., 1974.
\item" [2]" Treves F., Introduction to pseudo-differential and Fourier
integral operators. New York, 1980.
\item" [3]" Bogolyubov N.N., Mitropol'ski{\v\i} Yu.A., Asymptotical methods in
the theory of nonlinear oscillations [in Russian]. Moscow, Fizmatlit, 1958.
\item" [4]" Maslov V.P., Theory of perturbations and asymptotical methods
[in Russian]. Moscow, Izd-vo MGU, 1965.
\item" [5]" Maslov V.P., Fedoryuk M.V., Quasiclassical approximation for
equations of quantum mechanics [in Russian]. Moscow, Nauka, 1976.
\item" [6]" Karasev M.V., Maslov V.P., Nonlinear Poisson brackets. Geometry
and quantization. AMS, Providence, 1993.
\item" [7]" Kirillov A.A., Elements of representation theory. Springer, 1976.
\item" [8]" Barut A., Raczka R., Theory of group representations and
applications. PWN -- Polish Sci. Publ. Warszawa, 1977.
\item" [9]" Juriev D., Topics in hidden symmetries. V. E-print:
funct-an/9611003 (1996).
\item"[10]" Juriev D., On the infinite dimensional hidden symmetries. III.
$q_R$--conformal symmetries at $q_R\!\to\!\infty$ and Berezin-Karasev-Maslov
asymptotic quantization of $C^\infty(S^1)$. E-print: funct-an/9702002 (1997).
\item"[11]" Zhelobenko D.P., Representations of the reductive Lie algebras
[in Russian]. Moscow, Nauka, 1993.
\item"[12]" Dixmier J., Alg\`ebres enveloppantes. Gauthier-Villars, Paris,
1974.
\item"[13]" Juriev D.V., Complex projective geometry and quantum projective
field theory. Theor. Math.Phys. 101(3) (1994) 1387-1403.
\item"[14]" Berezin F.A., Quantization in complex symmetric spaces [in
Russian]. Izvestiya AN SSSR. Ser.matem. 39(2) (1975) 363-402.
\item"[15]" Klimek S., Lesniewski A., Quantum Riemann surfaces. I.
Commun.Math.Phys. 146 (1992) 103-122.
\item"[16]" Biedenharn L., Louck J., Angular momentum in quantum mechanics.
Theory and appli\-cations. Encycl.Math.Appl. V.8. Addison Wesley Publ. Comp.
1981.
\item"[17]" Fuchs D.B., Cohomology of infinite-dimensional Lie algebras [in
Russian]. Moscow, Nau\-ka, 1984.
\item"[18]" Juriev D., Noncommutative geometry, chiral anomaly in the quantum
projective ($\sltwo$--invariant) field theory and $\jltwo$--invariance.
J.Math.Phys. 33 (1992) 2819-2822, (E): 34 (1993) 1615.
\item"[19]" Juriev D., Remarks on nonassociative structures in quantum
projective field theory: the central extension $\jltwo$ of the double
$\sltwo+\sltwo$ of the simple Lie algebra $\sltwo$ and related topics.
Acta Appl.Math. 50 (1998) 191-196.
\item"[20]" Juriev D.V., Quantum projective field theory: quantum-field
analogs of the Euler-Arnol'd equations in projective $G$-hypermultiplets.
Theor.Math.Phys. 98(2) (1994) 147-161.
\item"[21]" Mal'tsev A.I., On a representation of nonassociative rings
[in Russian]. Uspekhi Matem. Nauk 7(1) (1952) 181-185 [reprinted in
``Selected papers. 1. Classical algebra'', Moscow, Nau\-ka, 1976, pp.328-331].
\item"[22]" Kobayashi Sh, Nomizu K., Foundations of differential geometry.
Interscience Publ., New York, 1963/69.
\item"[23]" Juriev D., Infinite dimensional geometry and quantum field theory
of strings. II. Infinite-dimensional noncommutative geo\-met\-ry of a
self-interacting string field. Russian J.Math. Phys. 4 (1996) 287-314.
\item"[24]" Juriev D., On the infinite dimensional hidden symmetries. I.
Infinite dimensional geo\-met\-ry of $q_R$--conformal symmetries.
E-print: funct-an/9612004 (1996).
\item"[25]" Juriev D., On the infinite dimensional hidden symmetries. II.
$q_R$--conformal modular functor. E-print: funct-an/9701009 (1997).
\item"[26]" Kirillov A.A., Introduction into the representation theory and
noncommutative harmonic analysis [in Russian]. Current Math. Problems.
Basic Directions. V.22, Moscow, VINITI, 1987, pp.5-162.
\item"[27]" Dubnov V.L., Maslov V.P., Naza{\v\i}kinski{\v\i} V.E.,
The complex Lagrangian germ and the canonical operator. Russian J.Math.Phys.
3 (1995) 141-180.
\item"[28]" Juriev D.V., Droems: experimental mathematics, informatics
and infinite-dimensional geometry [in Russian]. Report RCMPI-96/05$^+$
(1996).
\endroster

\newpage\tenpoint
\centerline{\bf APPENDIX: THE ORIGINAL RUSSIAN VERSION OF ARTICLE}
\ \newline
\cyr UDK: {\rm 517.986.5+512.554.32+512.813.4}\newline
\ \newline
\ \newline
\centerline{\cyb PRIBLIZHENNYE PREDSTAVLENIYA I ALGEBRA VIRASORO}
\ \newline
\ \newline
\centerline{\cyb D.V.Yurp1ev}
\ \newline\eightpoint
\centerline{\cyre Tsentr matematicheskoe0 fiziki i informatiki ``Talassa
E1teriya'',}
\centerline{\cyre ul.Miklukho-Maklaya 20-180, Moskva 117437 Rossiya.}
\centerline{E-mail: denis\@juriev.msk.ru}
\ \tenpoint\newline
\centerline{{\cyr Fevralp1, 1997}\qquad\qquad\qquad math.RT/9805001}
\ \newline
\ \newline

\cyr
Dannaya rabota, adresovannaya, v osnovnom, spetsialistam v teorii
pred\-stav\-le\-nie0 konechnomernyh reduktivnyh i beskonechnomernyh
$\Bbb Z$--gra\-du\-i\-ro\-van\-nyh algebr Li, teorii priblizhenie0 v
funktsionalp1nom analize i asim\-p\-to\-ti\-ches\-kim metodam v ischislenii
operatorov, a takzhe, otchasti, v ma\-te\-ma\-ti\-ches\-koe0 fizike
(algebraicheskoe0 kvantovoe0 teorii polya), po\-svya\-shche\-na izucheniyu
algebraicheskoe0 i analiticheskoe0 struktury bes\-ko\-nech\-no\-mer\-nyh
skrytyh simmetrie0 v teorii predstavlenie0 konechnomernyh re\-duk\-tiv\-nyh
algebr Li. Bolee tochno, predmetom raboty yavlyaet\-sya re\-a\-li\-za\-tsiya
beskonechnomernoe0 $\Bbb Z$--graduirovannoe0 algebry Vitta loranovskih
po\-li\-no\-mi\-alp1\-nyh vektornyh polee0 na okruzhnosti i ee odnomernogo
ne\-tri\-vi\-alp1\-no\-go tsentralp1nogo rasshireniya (algebry Virasoro)
skrytymi sim\-met\-ri\-ya\-mi v modulyah Verma nad algebroe0 Li $\sltwo$
(t.n. $q_R$--kon\-for\-m\-ny\-mi simmetriyami). Pri e1tom, ukazannye
beskonechnomernye algebry Li re\-a\-li\-zu\-yut\-sya $q_R$--konformnymi
simmetriyami ne tochno, a priblizhenno. V sta\-tp1e rassmatrivayut\-sya
dva tipa priblizhennyh predstavlenie0: pred\-s\-tav\-le\-niya po modulyu
nekotorogo klassa $\frak S$ operatorov (operatorov Gilp1\-ber\-ta-Shmidta,
kompaktnyh operatorov ili operatorov konechnogo ranga) i asimptoticheskie
predstavleniya ``$\mdohbar$'', gde $\hbar$ -- nekotorye0 pa\-ra\-metr (v
dannom sluchae vozmozhny razlichnye opredeleniya operatornogo
``$O$--bolp1\-sho\-go'' po parametru $\hbar$). Priblizhennye predstavleniya
pervogo ti\-pa estestvenno traktovatp1 v kontekste psevdodifferentsialp1nogo
is\-chis\-le\-niya [1,2], v to vremya kak vtorye -- asimptoticheskih metodov
[3-5]. Ot\-me\-tim, chto asimptoticheskie
predstavleniya ``$\mdohbar$'' issledovalisp1 v ramkah formalizma
asimptoticheskogo kvantovaniya Karaseva-Maslova [6], a predstavleniya po
modulyu $\frak S$ yavlyayut\-sya v opredelennom smysle obob\-shche\-ni\-ya\-mi
proektivnyh predstavlenie0 [7,8].

Rezulp1taty raboty, nekotorye iz kotoryh byli anonsirovany v
e1lek\-t\-ron\-nyh preprintah e1lektronnogo arhiva po funktsionalp1nomu
analizu Natsionalp1nyh laboratorie0 SShA v Los Alamose [9,10], nosyat skoree
harakter nablyudenie0, nezheli teorem, trebuyushchih trudoemkih i ne yasnyh
iz ih formulirovok dokazatelp1stv. Dokazatelp1stva predstavlennyh v sta\-tp1e
faktov po suti dela odnohodovy i predstavlyayut soboe0 gromozdkie
vy\-chis\-li\-telp1\-nye uprazhneniya (gde e1to vozmozhno, opushchennye),
opi\-ra\-yu\-shchi\-e\-sya na sovokupnostp1 poluchennyh avtorom ranee
yavnyh formul; takim ob\-ra\-zom, smyslovaya nagruzka prihodit\-sya, v
osnovnom, na sami for\-mu\-li\-rov\-ki re\-zulp1\-ta\-tov.

\head\cyb\S1. Predvaritelp1nye opredeleniya\endhead

\subhead\cyb 1.1. Algebra Li $\sltwo$ i moduli Verma nad nee0\endsubhead
\cyr Algebroe0 Li $\sltwo$ nazyvaet\-sya trehmernoe prostranstvo
kompleksnyh matrits $2\times2$ so sledom nulp1, snabzhennoe standartnym
kommutatorom $[X,Y]=XY-YX$, gde um\-no\-zhe\-nie v pravoe0 chasti --
standartnoe umnozhenie matrits. V bazise
$$ l_{-1}=\left(\matrix0&1\\0&0\endmatrix\right),\quad
l_0=\left(\matrix\tfrac12&0\\0&-\tfrac12\endmatrix\right),\quad
l_1=\left(\matrix0&0\\-1&0\endmatrix\right) $$
kommutatsionnye sootnosheniya imeyut vid: $[l_i,l_j]=(i-j)l_{i+j}$
($i,j=-1,0,1$). Algebra Li $\sltwo$ $\Bbb Z$--graduirovana:
$\deg(l_i)=-\ad(l_0)l_i=i$, gde $\ad(X)$ -- operator prisoedinennogo
dee0stviya v algebre Li: $\ad(X)Y=[X,Y]$. Kak sledstvie,
$\Bbb Z$--graduirovannye moduli nad $\sltwo$ yavlyayut\-sya
$l_0$--di\-a\-go\-nalp1\-ny\-mi. Vektor $v$ v $\Bbb Z$--graduirovannom
module nad algebroe0 Li $\sltwo$ nazyvaet\-sya e1kstremalp1nym, esli $l_1v=0$
i linee0naya obolochka vektorov $l_{-1}^nv$ ($n\in\Bbb Z_+$) sovpadaet s
samim modulem (t.e. $v$ -- tsiklicheskie0 vektor). $\Bbb Z$--graduirovannye0
modulp1 v kotorom sushchestvuet e1kstremalp1nye0 vektor (v e1tom sluchae on
opredelen s tochnostp1yu do umnozheniya na konstantu) nazyvaet\-sya
e1kstremalp1nym [11]. E1kstremalp1nye0 modulp1 nazyvaet\-sya mo\-du\-lem Verma
[12], esli dee0stvie $l_{-1}$ v nem svobodno, t.e. vektora $l_{-1}^nv$
linee0no nezavisimy. V sluchae algebry Li $\sltwo$ moduli Verma sutp1 v
tochnosti beskonechnomernye e1kstremalp1nye moduli. E1kstremalp1nym vesom
modulya Verma nazyvaet\-sya chislo $h$, opredelyaemoe iz ravenstva $l_0v=hv$,
gde $v$ -- e1kstremalp1nye0 vektor. Moduli Verma opredeleny dlya vseh
kompleksnyh chisel $h$ i poparno neizomorfny. V dalp1nee0shem my budem
rassmatrivatp1 moduli Verma tolp1ko s veshchestvennymi e1kstremalp1nymi
vesami.

Modulp1 Verma $V_h$ nad algebroe0 Li $\sltwo$ s e1kstremalp1nym vesom $h$
mozhet bytp1 realizovan v prostranstve $\Bbb C[z]$ mnogochlenov odnogo
kom\-p\-lek\-s\-no\-go peremennogo $z$. Formuly dlya generatorov algebry Li
$\sltwo$ imeyut vid:
$$l_{-1}=z,\quad l_0=z\partial_z+h,\quad l_1=z\partial_z^2+2h\partial_z,$$
zdesp1 $\partial_z=\frac{d}{dz}$.

Modulp1 Verma nevyrozhden (t.e. ne soderzhit sobstvennogo podmodulya) pri
$h\ne-\tfrac{n}2$ ($n\in\Bbb Z_+$). Modulp1 Verma $V_h$ nazyvaet\-sya
unitarizuemym (ili e1rmitovym), esli v nem zadana struktura predgilp1bertova
prostranstva takaya, chto $l_i^*=l_{-i}$. Popolnenie unitarizuemogo modulya
Verma budet oboznachatp1sya $V^{\Hilb}_h$. Algebra Li $\sltwo$ dee0stvuet
v $V^{\Hilb}_h$ ne\-og\-ra\-ni\-chen\-ny\-mi operatorami. Polezno takzhe
rassmatrivatp1 formalp1nye moduli Verma $V^{\form}_h$, kotorye realizuyut\-sya
v prostranstve $\Bbb C[[z]]$ formalp1nyh ste\-pen\-nyh ryadov odnogo
kompleksnogo peremennogo $z$, a formuly dlya ge\-ne\-ra\-to\-rov algebry Li
$\sltwo$ sovpadayut s privedennymi vyshe. Pri e1tom, $V_h\subseteq
V^{\Hilb}_h\subseteq V^{\form}_h$, i moduli $V_h$, $V^{\Hilb}_h$,
$V^{\form}_h$ obrazuyut troe0ku Gelp1fanda ili osnashchenie gilp1bertova
prostranstva $V^{\Hilb}_h$ po Diraku. Dee0stvie veshchestvennoe0 formy
algebry Li $\sltwo$, porozhdennoe0 antie1rmitovymi operatorami $il_0$,
$l_1-l_{-1}$ i $i(l_1+l_{-1})$ v gilp1bertovom
prostranstve $V^{\Hilb}_h$ neogranichennymi operatorami,
e1ks\-po\-nen\-tsi\-ru\-et\-sya do unitarnogo pred\-s\-tav\-le\-niya
sootvet\-stvuyushchee0 odnosvyaznoe0 gruppy Li.

V neunitarizuemom module Verma nad algebroe0 Li $\sltwo$ imeet\-sya
edinstvennaya (s tochnostp1yu do mnozhitelya) indefinitnaya
po\-lu\-to\-ra\-li\-nee0\-naya forma $(\cdot,\cdot)$ takaya, chto
$(l_iv_1,v_2)=(v_1,l_{-i}v_2)$ dlya lyubyh dvuh vektorov $v_1$ i $v_2$
iz modulya Verma. Esli e1ta polutoralinee0naya forma nevyrozhdena ( a v
e1tom sluchae i sam modulp1 nevyrozhden), to ona imeet signaturu vida
$(n,\infty)$, gde $n$ -- konechno, a sledovatelp1no, opredeleno popolnenie
takogo modulya Verma po Pontryaginu [10]. Sootvet\-stvuyushchie0 modulp1,
v kotorom algebra Li $\sltwo$ dee0stvuet neogranichennymi operatorami,
bu\-det oboznachatp1sya $V^{\Pontr}_h$. Imeet mesto tsepochka vklyuchenie0:
$V_h\subseteq V^{\Pontr}_h\subseteq V^{\form}_h$. Dee0stvie0
veshchestvennoe0 formy algebry Li $\sltwo$, porozhdennoe0 antie1rmitovymi
(otnositelp1no nevyrozhdennoe0 indefinitnoe0 po\-lu\-to\-ra\-li\-nee0\-noe0
formy $(\cdot,\cdot)$) operatorami $il_0$, $l_1-l_{-1}$ i $i(l_1+l_{-1})$
v prostranstve Pontryagina $V^{\Pontr}_h$ neogranichennymi
ope\-ra\-to\-ra\-mi, e1ksponentsiruet\-sya do psevdounitarnogo predstavleniya
sootvet\-stvuyushchee0 odnosvyaznoe0 gruppy Li.

\subhead\cyb 1.2. Skrytye simmetrii v modulyah Verma nad algebroe0 Li
$\sltwo$: $C^*$--algebra Lobachevskogo-Berezina i $q_R$--konformnye
simmetrii\endsubhead

\proclaim{\cyb Predlozhenie 1 [13]}\cyi V nevyrozhdennom module Verma
$V_h$ nad algebroe0 Li $\sltwo$ odnoznachno opredeleny operatory $D$ i
$F$ takie, chto
$$\aligned [D,l_{-1}]=1,\quad [D,l_0]=D,\quad [D,l_1]=D^2,\\
[l_{-1},F]=1,\quad [l_0,F]=F,\quad [l_1,F]=F^2.\endaligned$$
Esli moduli Verma realizovany v prostranstve $\Bbb C[z]$ mnogochlenov odnoe0
kompleksnoe0 peremennoe0 $z$, to
$$D=\partial_z,\qquad F=z\tfrac1{\xi+2h},$$
gde $\xi=z\partial_z$. Operatory $F$ i $D$ udovletvoryayut sootnosheniyam:
$$[FD,DF]=0,\qquad [D,F]=q_R(1-DF)(1-FD),$$
gde $q_R=\tfrac1{2h-1}$. V unitarizuemom module Verma (pri $q_R\!\ne\!0$)
operatory $F$ i $D$ ogranicheny i $F^*=D$, $D^*=F$.
\endproclaim

\cyr
Algebra, porozhdennaya peremennymi $t$ i $t^*$, s sootnosheniyami
$[tt^*,t^*t]=0$ i $[t,t^*]=q_R(1-tt^*)(1-t^*t)$, buduchi kvantovaniem po
Berezinu ploskosti Lobachevskogo, realizovannoe0 v edinichnom kompleksnom
diske (re\-a\-li\-za\-tsiya Puankare) [14], nazyvaet\-sya algebroe0
Lobachevskogo-Berezina. Pred\-lo\-zhe\-nie 1 pozvolyaet rassmatrivatp1
algebru Lobachevskogo-Berezina kak $C^*$--algebru. $C^*$--algebra
Lobachevskogo-Berezina byla nedavno pe\-re\-otk\-ry\-ta S.Klimekom i
A.Lesnievskim [15].

\proclaim{\cyb Predlozhenie 2 [13]}\cyi V nevyrozhdennom module Verma
$V_h$ nad algebroe0 Li $\sltwo$ odnoznachno opredeleny operatory $L_n$
($n\in\Bbb Z$) takie, chto
$$[l_i,L_n]=(i-n)L_{i+n}\qquad (i=1,2,1;\ n\in\Bbb Z),$$
pri e1tom, $L_i=l_i$ ($i=-1,0,1$). Esli moduli Verma realizovany v
prost\-ran\-st\-ve $\Bbb C[z]$ mnogochlenov odnoe0 kompleksnoe0 peremennoe0
$z$, to
$$L_k=(xi+(k+1)h)\partial_z^k\ (k\ge0),\quad
L_{-k}=z^k\tfrac{\xi+(k+1)h}{(\xi+2h)\ldots(\xi+2h+k-1)}\ (k\ge1),$$
gde $\xi=z\partial_z$. Operatory $L_n$ udovletvoryayut sootnosheniyam:
$$[L_n,L_m]=(n-m)L_{n+m},\text{\cyr\ esli}\ n,m\ge-1\text{\cyr\ ili}\ n,m\le1.$$
V unitarizuemom module Verma operatory $L_n$ neogranicheny i $L_i^*=L_{-1}$.
\endproclaim

\cyr
Operatory $L_n$ nazyvayut\-sya $q_R$--konformnymi simmetriyami. Oni {\cyi
sim\-vo\-li\-ches\-ki\/} mogut bytp1 predstavleny v vide:
$$L_n=D^{nh}L_0D^{n(1-h)},\qquad L_{-n}=F^{n(1-h)}L_0F^{nh}.$$
Dlya pridaniya smysla simvolicheskoe0 zapisi sleduet vospolp1zovatp1sya
ob\-shchi\-mi perestanovochnymi sootnosheniyami
$$[L_n,f(D)]=(-D)^{n+1}f'(D)\ (n\ge-1),\qquad
[L_{-n},f(F)]=F^{n+1}f'(D)\ (n\ge-1)$$
pri $n=0$.

Sootnosheniya kommutatsii operatorov $D$, $F$ i generatorov
$q_R$--kon\-form\-nyh simmetrie0 s generatorami algebry Li $\sltwo$ oznachayut,
chto se\-mee0\-s\-t\-va $J_k$ i $L_k$ ($k\in\Bbb Z$), gde $J_i=D^i$,
$J_{-i}=F^i$ ($i\in\Bbb Z_+$), yavlyayut\-sya semee0stvami tenzornyh
operatorov [8,16] dlya algebry Li $\sltwo$.

S semee0stvami tenzornyh operatorov $J_k$ i $L_k$ mozhno assotsiirovatp1 ih
proizvodyashchie funktsii -- operatornye polya, t.e. formalp1nye ryady
Lorana ot odnoe0 kompleksnoe0 peremennoe0 $u$ s operatornymi
ko\-e1f\-fi\-tsi\-en\-ta\-mi:
$$J(u)=\sum_{i\in\Bbb Z}J_i(-u)^{-1-i},\quad
T(u)=\sum_{i\in\Bbb Z}L_i(-u)^{-2-i}.$$
Pri kazhdom znachenii $u$ operatornoe pole zadaet otobrazhenie $V_h$ v
$V^{\form}_h$. Tot fakt, chto $J_k$ i $T_k$ ($k\in\Bbb Z$) obrazuyut
semee0stva tenzornyh operatorov, na yazyke operatornyh polee0 oznachaet, chto
$$[l_i,J(u)]=(-u)^iJ(u)-(-u)^{i+1}J'(u),\quad
[l_i,T(u)]=2(-u)^iT(u)-(-u)^{i+1}T'(u).$$
Operatornye polya $V(u)$, preobrazuyushchiesya kak $s$--differentsialy pod
dee0\-s\-t\-vi\-em algebry Li $\sltwo$ (t.e. udovletvoryayushchie tozhdestvu
$[l_i,V(u)]=s(-u)^iV(u)-(-u)^{i+1}V'(u)$), nazyvayut\-sya
$\sltwo$--pervichnymi spi\-na $s$; takim obrazom, operatornye polya $J(u)$
i $T(u)$ yavlyayut\-sya $\sltwo$--pervichnymi spi\-na 1 i 2, sootvet\-stvenno.
Operatornye polya $J(u)$ i $T(u)$, a takzhe ih svoe0stva podrobno izuchalisp1
v [13].

\subhead\cyb 1.3. Beskonechnomernye $\Bbb Z$--graduirovannye algebry Li:
algebra Vit\-ta $\frak w^{\Bbb C}$ loranovskih polinomialp1nyh vektornyh
polee0 na okruzhnosti i algebra Virasoro $\vir^{\Bbb C}$, ee odnomernoe
netrivialp1noe tsentralp1noe rasshirenie
\endsubhead
\cyr Algebra Li $\Vect(S^1)$ realizuet\-sya v prostranstve
$C^\infty$--glad\-kih vektornyh polee0 $v(t)\partial_t$ na okruzhnosti
$S^1\simeq\Bbb R/2\pi\Bbb Z$ s kommutatorom
$$[v_1(t)\partial_t,v_2(t)\partial_t]=
(v_1(t)v_2'(t)-v_1'(t)v_2(t))\partial_t.$$
V bazise $s_n=\sin(nt)\partial_t$, $c_n=\cos(nt)\partial_t$, $h=\partial_t$
kommutatsionnye sootnosheniya imeyut vid:
$$\aligned
[s_n,s_m]&=\tfrac12((m-n)s_{n+m}+\sgn(n-m)(n+m)s_{|n-m|}),\\
[c_n,c_m]&=\tfrac12((n-m)s_{n+m}+\sgn(n-m)(n+m)s_{|n-m|}),\\
[s_n,c_m]&=\tfrac12((m-n)c_{n+m}-(m+n)c_{|n-m|})-n\delta_{nm}h,\\
[h,s_n]&=nc_n,\quad [h,c_n]=-ns_n.
\endaligned$$
Oboznachim $\Vect^{\Bbb C}(S^1)$ kompleksifikatsiyu algebry Li $\Vect(S^1)$.
V bazise $e_n=ie^{ikt}\partial_t$ kommutatsionnye sootnosheniya v algebre
Li $\Vect^{\Bbb C}(S^1)$ imeyut vid:
$$[e_j,e_k]=(j-k)e_{j+k}.$$

Udobno takzhe rassmatrivatp1 vlozhenie okruzhnosti $S^1$ v kompleksnuyu
ploskostp1 $\Bbb C$ s koordinatoe0 $z$, pri e1tom na okruzhnosti $z=e^{it}$,
a e1lementy bazisa $e_k$ ($k\in\Bbb Z$) predstavlyayut\-sya loranovskimi
polinomialp1nymi vek\-tor\-ny\-mi polyami: $e_k=z^{k+1}\partial_z$.
$\Bbb Z$--graduirovannaya algebra Li, po\-rozh\-den\-naya loranovskimi
polinomialp1nymi vektornymi polyami (t.e. konechnymi li\-nee0\-ny\-mi
kombinatsiyami e1lementov bazisa $e_k$), nazyvaet\-sya algebroe0 Vit\-ta
i oboznachaet\-sya $\frak w^{\Bbb C}$. Algebra Vitta $\frak w^{\Bbb C}$
yavlyaet\-sya kompleksifikatsiee0 podalgebry $\frak w$ algebry $\Vect(S^1)$,
porozhdennoe0 trigonometricheskimi po\-li\-no\-mi\-alp1\-ny\-mi vektornymi
polyami na okruzhnosti $S^1$, t.e. konechnymi linee0nymi kombinatsiyami
e1lementov bazisa $s_n$, $c_n$ i $h$.

Algebra Li $\Vect(S^1)$ dopuskaet netrivialp1noe odnomernoe tsentralp1noe
rasshirenie, zadavaemoe 2-kotsiklom Gelp1fanda-Fuksa [17]:
$$c(v_1(t)\partial_t,v_2(t)\partial_t)=
\int_0^{2\pi}(v_1'(t)v_2''(t)-v_2'(t)v_1'(t))\,dt.$$
Dannoe rasshirenie, buduchi prodolzhennym na kompleksifikatsiyu
$\Vect^{\Bbb C}(S^1)$ algebry Li $\Vect(S^1)$ i ogranichennym na podalgebru
$\frak w^{\Bbb C}$, zadaet od\-no\-mer\-noe tsentralp1noe rasshirenie algebry
Vitta, nazyvaemoe algebroe0 Virasoro i oboznachaemoe $\vir^{\Bbb C}$.
Algebra Virasoro porozhdaet\-sya ge\-ne\-ra\-to\-ra\-mi $e_k$ ($k\in\Bbb Z$)
i tsentralp1nym e1lementom $c$ s kommutatsionnymi sootnosheniyami:
$$[e_j,e_k]=(j-k)e_{j+k}+\frac{j^3-j}{12}c$$
i yavlyaet\-sya kompleksifikatsiee0 tsentralp1nogo rasshireniya $\vir$
algebry Li $\frak w$. V neprivodimom predstavlenii tsentralp1nye0 e1lement
$c$ algebry Vi\-ra\-so\-ro perehodit v skalyarnye0 operator, koe1ffitsient
proportsionalp1nosti kotorogo edinichnomu operatoru nazyvaet\-sya
{\cyi tsentralp1nym zaryadom}.

\head\cyb\S2. $\HS$--proektivnye predstavleniya algebry Vitta v modulyah
Verma nad algebroe0 Li $\sltwo$\endhead

\subhead\cyb 2.1. $\frak A$--proektivnye predstavleniya [9]\endsubhead

\definition{\cyb Opredelenie 1A} \cyr Pustp1 $\frak A$ -- proizvolp1naya
assotsiativnaya algebra, pred\-s\-tav\-len\-naya linee0nymi operatorami
v linee0nom prostranstve $H$ i $\frak g$ -- algebra Li. Linee0noe
otobrazhenie $T:\frak g\mapsto\End(H)$ nazyvaet\-sya {\cyi
$\frak A$--pro\-ek\-tiv\-nym predstavleniem}, esli dlya proizvolp1nyh $X$
i $Y$ iz $\frak g$ sushchestvuet e1lement algebry $\frak A$, predstavlennye0
operatorom $A_{XY}$, takoe0, chto
$$[T(X),T(Y)]-T([X,Y])=A_{XY}.$$
\enddefinition

Esli $H$ beskonechnomerno, to predstavlenie mozhet osushchestvlyatp1sya
ne\-og\-ra\-ni\-chen\-ny\-mi operatorami s podhodyashchimi oblastyami
opredeleniya.

\remark{\cyi Zamechanie 1} \cyr Opredelenie 1A mozhet bytp1 obobshcheno na
proizvolp1nye an\-ti\-kom\-mu\-ta\-tiv\-nye algebry. V e1toe0 situatsii ono
tesno svyazano s kon\-s\-t\-ruk\-tsi\-ya\-mi predstavlenie0 antikommutativnyh
algebr $\jltwo$ i $\slstwo$ v [18,19;20:\S2]. V tselom, ego sleduet
rassmatrivatp1 v kontekste staryh idee0 A.I.Malp1tseva o predstavleniyah
proizvolp1nyh neassotsiativnyh al\-gebr [21]. Obshchie antikommutativnye
algebry i konstruktsii ih $\frak A$--pro\-ek\-tiv\-nyh predstavlenie0
interesny s tochki zreniya kvazihopfovyh algebr, neassotsiativnyh kak
koalgebry, yakobievyh i koyakobievyh kvazibialgebr i svyazannyh s nimi
struktur (sm.ssylki v [9]).
\endremark

\remark{\cyi Primer} \cyr Esli $(\frak g,\frak h)$ -- reduktivnaya para, to
lyuboe predstavlenie algebry Li $\frak g$ yavlyaet\-sya
$\Cal U(\frak h)$--proektivnym predstavleniem binarnoe0
an\-ti\-kom\-mu\-ta\-tiv\-noe0 algebry $\frak p$ ($\frak g=\frak h\oplus\frak
p$, binarnaya operatsiya v $\frak p$ imeet standartnye0 vid:
$[X,Y]_{\frak p}=\pi([X,Y])$, gde $[\cdot,\cdot]$ -- kommutator v $\frak g$,
a $\pi$ -- proektor $\frak g$ na $\frak p$ vdolp1 $\frak h$ [22]).
\endremark

\remark{\cyi Zamechanie 2} \cyr Standartnoe proektivnoe predstavlenie
yavlyaet\-sya chastnym sluchaem opredeleniya, esli algebra $\frak A$
odnomerna i dee0stvuet v $H$ ska\-lyar\-ny\-mi matritsami.
\endremark

\remark{\cyi Zamechanie 3} \cyr Esli $H$ -- gilp1bertovo (ili
predgilp1bertovo) prostranstvo, to v kachestve $\frak A$ mozhno
rassmatrivatp1 algebru $\HS$ vseh operatorov Gilp1\-ber\-ta-Shmidta. Mozhno
takzhe rassmatrivatp1 algebry $\bnd$, $\K$, $\TC$ i $\FR$ ogranichennyh,
kompaktnyh, yadernyh operatorov i operatorov konechnogo ranga.
\endremark

\definition{\cyb Opredelenie 1B} \cyr Pustp1 $\frak A$ -- assotsiativnaya
algebra s involyutsiee0 $*$, simmetrichno predstavlennaya v gilp1bertovom
prostranstve $H$. Esli $\frak g$ -- algebra Li s involyutsiee0 $*$, to ee
$\frak A$--proektivnoe predstavlenie $T$ v prostranstve $H$ nazyvaet\-sya
{\cyi simmetrichnym}, esli dlya vseh e1lementov $a$ iz $\frak g$
$T(a^*)=T^*(a)$. Pustp1 $\frak g$ -- $\Bbb Z$--graduirovannaya algebra Li
($\frak g=\oplus_{n\in\Bbb Z}\frak g_n$) s involyutsiee0 $*$ takaya, chto
$\frak g^*_n=\frak g_{-n}$ i involyutsiya tozhdestvenna na podalgebre
$\frak g_0$. Prodolzhim $\Bbb Z$--graduirovku i involyutsiyu $*$ s algebry
Li $\frak g$ na tenzornuyu algebru $\bold T^{\cdot}(\frak g)$. Simmetrichnoe
$\frak A$--proektivnoe predstavlenie $\frak g$ nazyvaet\-sya {\cyi
absolyutno simmetrichnym}, esli dlya lyubogo e1lementa $a$ al\-geb\-ry $\bold
T^{\cdot}(\frak g)$ takogo, chto $\deg(a)=0$, vypolneno ravenstvo
$T(a)=T^*(a)$ (zdesp1 predstavlenie $T$ algebry Li $\frak g$ v $H$
prodolzheno do otobrazheniya $\bold T^{\cdot}(\frak g)$ v $\End(H)$.
\enddefinition

\definition{\cyb Opredelenie 1V} \cyr $\frak A$--proektivnoe predstavlenie
$T$ algebry Li $\frak g$ v li\-nee0\-nom prostranstve $H$ nazyvaet\-sya {\cyi
pochti absolyutno zamknutym}, esli dlya lyubogo naturalp1nogo $n$ i
proizvolp1nyh e1lementov $X_0, X_1, X_2,\ldots X_{n+1}$ algebry Li $\frak
g$ sushchestvuet e1lement $\varphi(X_0,X_1,X_2,\ldots X_{n+1})$ algebry
$\frak g$ takoe0, chto
$$[\ \ldots[[T(X_0),T(X_1)],T(X_2)],\,\ldots,T(X_{n+1})]\!\equiv\!
T(\varphi(X_0,X_1,X_2,\ldots X_{n+1}))\!\!\!\!\pmod{\frak A},$$
zdesp1 $\frak A$ rassmatrivaet\-sya predstavlennoe0 v $\End(H)$. Pochti
absolyutno zamknutoe $\frak A$--proektivnoe predstavlenie $T$ algebry Li
$\frak g$ v linee0nom pros\-t\-ran\-st\-ve $H$ nazyvaet\-sya {\cyi absolyutno
zamknutym}, esli $\varphi(\cdot,\ldots,\cdot)\equiv 0$.
\enddefinition

\cyr Otobrazheniya
$(X_0,X_1,X_2,\ldots X_{n+1})\mapsto\varphi(X_0,X_1,X_2,\ldots X_{n+1})$,
as\-so\-tsi\-i\-ro\-van\-nye s proizvolp1nym pochti absolyutno zamknutym
$\frak A$--proektivnym pred\-s\-tav\-le\-ni\-em algebry Li $\frak g$,
opredelyayut vysshie skobki v algebre Li $\frak g$. Obp2ekty s vysshimi
skobkami sistematicheski poyavlyayut\-sya vo mnogih ob\-las\-tyah matematiki i
matematicheskoe0 fiziki (sm.napr.knigu [6] ili sta\-tp1yu [23] sredi prochih
i mnogochislennye ssylki v nih).

\subhead\cyb 2.2. $\HS$--proektivnye predstavleniya algebry Vitta
$q_R$--konformnymi simmetriyami v unitarizuemyh modulyah Verma $V_h$ nad
algebroe0 Li $\sltwo$ [9]
\endsubhead
\cyr Otmetim, chto algebra Vitta $\frak w^{\Bbb C}$ dopuskaet estestvennuyu
in\-vo\-lyu\-tsiyu $*$.

\proclaim{\cyb Teorema 1A}\cyi Generatory $L_k$ ($k\in\Bbb Z$)
$q_R$--konformnyh simmetrie0 v uni\-ta\-ri\-zu\-e\-mom module Verma $V_h$
nad algebroe0 Li $\sltwo$ osushchestvlyayut ab\-so\-lyut\-no simmetrichnoe
$\HS$--proektivnoe predstavlenie algebry Vitta $\frak w^{\Bbb C}$. Dobavlenie
k tenzornym operatoram $L_k$ tenzornyh operatorov $J_k$ ($k\in\Bbb Z$)
privodit k absolyutno simmetrichnomu $\HS$--proektivnomu predstavleniyu
polupryamoe0 summy algebry Vitta i bes\-ko\-nech\-no\-mer\-noe0 $\Bbb
Z$--graduirovannoe0 algebry Gee0zenberga (odnomernogo tsentralp1nogo
rasshireniya bes\-ko\-nech\-no\-mer\-noe0 $\Bbb Z$--graduirovannoe0
abelevoe0 algebry Li $\Bbb C[z,z^{-1}]$ mnogochlenov Lorana).
\endproclaim

\cyr Utverzhdenie teoremy sleduet iz yavnyh formul dlya generatorov
$q_R$--konformnyh simmetrie0 i tenzornyh operatorov $J_k$. Proverka na
pri\-nad\-lezh\-nostp1 klassu $\HS$ operatorov Gilp1berta-Shmidta dlya
operatorov fiksirovannoe0 stepeni otnositelp1no graduirovki v $\Bbb
Z$--graduirovannom prostranstve mnogochlenov $\Bbb C[z]$, snabzhennom
nekotorym skalyarnym pro\-iz\-ve\-de\-ni\-em, otnositelp1no kotorogo
odnomernye vesovye prostranstva or\-to\-go\-nalp1\-ny, kakovym i
yavlyaet\-sya unitarizuemye0 modulp1 Verma nad algebroe0 Li $\sltwo$,
ne predstavlyaet nikakih problem.

\remark{\cyi Zamechanie 4} \cyr $\HS$--proektivnye predstavleniya algebry
Vitta v uni\-ta\-ri\-zu\-e\-myh modulyah Verma nad $\sltwo$ absolyutno
zamknuty.
\endremark

\remark{\cyi Zamechanie 5} \cyr Teorema 1A perenosit\-sya na psevdounitarnye0
sluchae0 s za\-me\-noe0 klassa $\HS$ operatorov Gilp1berta-Shmidta na klass
$\K$ kompaktnyh operatorov.
\endremark

Rezulp1taty, kasayushchiesya ``e1ksponentsirovannoe0'' versii teoremy 1A,
anonsirovany v e1lektronnyh preprintah [24,25].

\subhead\cyb 2.3. $\FR$--proektivnye predstavleniya algebry Vitta v modulyah
Verma nad algebroe0 Li $\sltwo$ s e1kstremalp1nym vesom $1$ i $\frac12$
\endsubhead
\cyr V nekotoryh chastnyh sluchayah utverzhdenie teoremy 1A mozhet bytp1
usileno.

\proclaim{\cyb Teorema 1B}\cyi Pri $h=1$ ili $h=\frac12$ $\HS$--proektivnoe
predstavlenie algebry Vitta $q_R$--konformnymi simmetriyami v unitarizuemyh
modulyah Verma nad algebroe0 Li $\sltwo$ yavlyaet\-sya $\FR$--proektivnym.
Pri $h=\frac12$ $\FR$--pro\-ek\-tiv\-nym yavlyaet\-sya i $\HS$--proektivnoe
predstavlenie polupryamoe0 summy algebry Vitta i beskonechnomernoe0 algebry
Gee0zenberga.
\endproclaim

\cyr Utverzhdenie sleduet iz yavnyh formul dlya tenzornyh operatorov $L_k$
i $J_k$.

\remark{\cyi Zamechanie 6}\cyr Pri $h=0$ $\K$--proektivnoe predstavlenie
algebry Vitta v module Verma nad algebroe0 Li $\sltwo$ yavlyaet\-sya
$\FR$--proektivnym.
\endremark

\subhead\cyb 2.4. $\bnd$--proektivnye predstavleniya algebry Vitta\endsubhead
\cyr Otmetim, chto k chislu ``neudobstv'' mnogih $\frak A$--proektivnyh
predstavlenie0, v tom chisle esli $\frak A$ -- nekotorye0 klass operatorov
$\K$, $\HS$, $\TC$ ili $\FR$, otnosit\-sya, voobshche govorya,
nezamknutostp1 ih sovokupnosti otnositelp1no vzyatiya ten\-zor\-nyh
proizvedenie0. Dlya klassa zhe $\bnd$ ogranichennyh operatorov ten\-zor\-nye
proizvedeniya $\bnd$--proektivnyh predstavlenie0 opredeleny. Odnako,
vsyakoe $\frak S$--proektivnoe predstavlenie, gde $\frak S$ -- odin iz
upomyanutyh vyshe klassov, yavlyaet\-sya $\bnd$--proektivnym i kak takovoe,
inogda, netrivialp1nym (kogda is\-hodnye predstavleniya osushchestvlyalisp1
neogranichennymi ope\-ra\-to\-ra\-mi), chto pozvolyaet konstruirovatp1 ih
tenzornye proizvedeniya, ko\-to\-rye budut v e1tom sluchae netrivialp1nymi
$\bnd$--proektivnymi pred\-s\-tav\-le\-ni\-ya\-mi.

Nazovem $\frak S$--proektivnoe predstavlenie algebry Li $\frak g$ ($\frak S$
-- nekotorye0 operatornye0 klass) {\cyi $\frak S$--neprivodimym}, esli
operatory predstavleniya od\-no\-vre\-men\-no ne mogut bytp1 privedeny putem
dobavleniya k nim e1lementov iz $\frak S$ k vidu, v kotorom oni vse
obladali by obshchim dlya nih sobstvennym invariantnym podprostranstvom.

\proclaim{\cyb Teorema 2}\cyi Pustp1 $T_h$ -- $\bnd$--proektivnoe
predstavlenie algebry Vitta $\frak w^{\Bbb C}$ $q_R$--konformnymi
simmetriyami v module Verma $V_h$ nad algebroe0 Li $\sltwo$, togda
$\bnd$--proektivnye predstavleniya $S^n(T_h)$ (zdesp1 $S^n$ oboznachaet
funktor vzyatiya $n$--oe0 simmetricheskoe0 stepeni, sr.[26])
$\bnd$--neprivodimy.
\endproclaim

\cyr Utverzhdenie teoremy sleduet iz yavnyh formul dlya generatorov
$q_R$--konformnyh simmetrie0.

\remark{\cyi Zamechanie 7}\cyr Razlozheniya tenzornyh proizvedenie0
$T_{h_1}\otimes\ldots\otimes T_{h_n}$ $\bnd$--pro\-ek\-tiv\-nyh predstavlenie0
$T_{h_i}$ na $\bnd$--neprivodimye komponenty v obshchem sluchae
neizvestny.
\endremark

\head\cyb\S3. Asimptoticheskie predstavleniya algebr Vitta i Virasoro v
modulyah Verma nad algebroe0 Li $\sltwo$\endhead

\subhead\cyb 3.1. Asimptotika $q_R$--konformnyh simmetrie0 v modulyah
Verma $V_h$ nad algebroe0 Li $\sltwo$ pri $q_R\to0$ ($h\to\frac12$) i
$q_R\to1$ ($h\to1$) i algebra Vitta $\frak w^{\Bbb C}$\endsubhead
\cyr Pustp1 $\frak A$ -- algebra konechnyh linee0nyh kombinatsie0
vyrazhenie0 vida $f(qp)p^n$ ($n\in\Bbb Z_+$) ili $q^nf(qp)$ ($n\in\Bbb Z_+$),
где $[p,q]=1$, a $f$ -- ratsionalp1nye funktsii, chp1i znamenateli ne
imeyut nulee0 v tselyh tochkah. Budem nazyvatp1 $\frak A$ {\cyi rasshirennoe0
algebroe0 Vee0lya}. Rasshirennaya algebra Vee0lya yavlyaet\-sya
topologicheskoe0 algebroe0 otnositelp1no estestvennoe0 s\-hodimosti.

\proclaim{\cyb Lemma}\cyi Rasshirennaya algebra Vee0lya $\frak A$ dopuskaet
tochnoe predstavlenie v module Verma $V_h$ $(h\not\in\Bbb Z/2)$:
$$p\to\partial_z,\qquad q\to z.$$
Pri e1tom, generatory $q_R$--konformnyh simmetrie0 $L_k$ $(k\in\Bbb Z)$,
ravno kak i generatory $D$, $F$ algebry Lobachevskogo-Berezina,
prinadlezhat obrazu ras\-shi\-ren\-noe0 algebry Vee0lya $\frak A$.
\endproclaim

\cyr Budem ponimatp1 v dalp1nee0shem ``$O$--bolp1shoe'' dlya asimptotik
ope\-ra\-to\-rov v modulyah Verma nad algebroe0 Li $\sltwo$ s nepolutselymi
e1ks\-t\-re\-malp1\-ny\-mi vesami, prinadlezhashchih rasshirennoe0 algebre
Vee0lya (po\-skolp1\-ku ee predstavleniya v modulyah Verma tochny mozhno
otozhdestvlyatp1 samu al\-geb\-ru s ee obrazom), v smysle topologii v
rasshirennoe0 algebre Vee0lya.

\proclaim{\cyb Teorema 3A}\cyi Generatory $q_R$--konformnyh simmetrie0 v
modulyah Verma $V_{\frac12+\hbar}$ $(0\!<\!\hbar\!<\!\frac12)$ nad algebroe0
Li $\sltwo$ obrazuyut asimptoticheskoe ``$\mdohbar$'' predstavlenie algebry
Vitta $\frak w^{\Bbb C}$.

Generatory $q_R$--konformnyh simmetrie0 v modulyah Verma $V_{1\pm\hbar}$
$(0\!<\!\hbar\!<\!\frac12)$ nad algebroe0 Li $\sltwo$ takzhe obrazuyut
asimptoticheskoe ``$\mdohbar$'' predstavlenie algebry Vitta $\frak w^{\Bbb C}$.
\endproclaim

\cyr Pervoe utverzhdenie teoremy, nemedlenno sleduyushchee iz yavnogo vida
generatorov $q_R$--konformnyh simmetrie0, bylo po suti dela dokazano v [10],
gde rassmatrivalasp1 sootvet\-stvuyushchaya ``e1ksponentsirovannaya''
ver\-siya. Vtoroe emu polnostp1yu podobno.

\remark{\cyi Zamechanie 8}\cyr Utverzhdenie teoremy perestaet bytp1 vernym,
elsi ``$O$--bolp1shoe'' rassmatrivaet\-sya v smysle (slaboe0) operatornoe0
s\-hodimosti v prostranstve operatorov v modulyah Verma $V_h$ nad algebroe0
Li $\sltwo$, otozhdestvlennyh s prostranstvom $\Bbb C[z]$ polinomov.
\endremark

\cyr Otvet na vopros o prichinah ``ras\-hozhdeniya'' mezhdu teoremoe0 3A
i za\-me\-cha\-ni\-em 8 zaklyuchaet\-sya v sleduyushchem: tochnye
predstavleniya rasshirennoe0 algebry Vee0lya ne tolp1ko ne sushchestvuyut
pri predelp1nom znachenii $\hbar\!=\!0$, no buduchi nepreryvnymi (v slaboe0
operatornoe0 topologii v prostranstve $V_h$) pri vseh znacheniyah parametra
$\hbar$ iz intervala $0\!<\!\hbar\!<\!\frac12$, ne yavlyayut\-sya {\cyi
ravnomerno\/} nepreryvnymi na e1tom intervale. Otmetim takzhe, chto
nepreryvnye predstavleniya rasshirennoe0 algebry Vee0lya v modulyah Ver\-ma
$V_h$ nad algebroe0 Li $\sltwo$ prodolzhayut\-sya do predstavlenie0
ne\-og\-ra\-ni\-chen\-ny\-mi operatorami v prostranstvah $V^{\Hilb}_h$
(ili $V^{\Pontr}_h$), i, kak sled\-s\-t\-vie, ne nepreryvnyh, chto, odnako,
tipichno dlya teorii predstavlenie0 algebr Li.

\subhead\cyb 3.2. Otsenka ``$\md\HS$'' vtorogo chlena asimptotiki
$q_R$--konformnyh simmetrie0 v modulyah Verma $V_h$ nad algebroe0 Li $\sltwo$
pri $q_R\!\to\!0$ ($h\!\to\!\frac12$) i pri $q_R\!\to\!1$ ($h\!\to\!1$) i
algebra Virasoro $\vir^{\Bbb C}$\endsubhead
\cyr Pri ``gib\-ri\-di\-za\-tsii'' priblizhenie0 ``$\md O(\hbar^n)$'' i
``$\md\frak S$'' nablyudaet\-sya interesnoe yavlenie ``nekommutiruemosti''
otsenok, primer chemu predostavlyayut asim\-p\-to\-ti\-ki $q_R$--konformnyh
simmetrie0.

\proclaim{\cyb Teorema 3B}\cyi Dlya generatorov $q_R$--konformnyh
simmetrie0 v modulyah Ver\-ma $V_{\frac12+\hbar}$ $(0\!<\!\hbar\!<\!\frac12)$
nad algebroe0 Li $\sltwo$ imeet\-sya sleduyushchaya asimptotika
``$\md O(\hbar^2)$'':
$$[L_i,L_j]=(i-j)L_{i+j}+\hbar\Cal A+O(\hbar^2),$$
gde
$$\Cal A\equiv\frac23(i^3-i)\delta_{ij} \mod{\HS}.$$
Takim obrazom, generatory $q_R$--konformnyh simmetrie0 obrazuyut
pri\-bli\-zhen\-noe predstavlenie algebry Virasoro $\vir^{\Bbb C}$ pri
$h\to\frac12$ v ukazannom smys\-le.

Dlya generatorov $q_R$--konformnyh simmetrie0 v modulyah Verma $V_{1\pm\hbar}$
$(0\!<\!\hbar\!<\!\frac12)$ nad algebroe0 Li $\sltwo$ takzhe imeet\-sya
asimptotika
``$\md O(\hbar^2)$'':
$$[L_i,L_j]=(i-j)L_{i+j}+\hbar\Cal A+O(\hbar^2),$$
gde
$$\Cal A\equiv\frac23(i^3-i)\delta_{ij} \mod{\HS}.$$
Takim obrazom, generatory $q_R$--konformnyh simmetrie0 obrazuyut
pri\-bli\-zhen\-noe predstavlenie algebry Virasoro $\vir^{\Bbb C}$ pri $h\to1$
v ukazannom smysle.
\endproclaim

\cyr V utverzhdenii teoremy otsenki ``$\md O(\hbar^2)$'' ponimayut\-sya
v smysle s\-hodimosti v rasshirennoe0 algebre Vee0lya $\frak A$, pri e1tom
$\HS$ oboznachaet sovokupnostp1 e1lementov e1toe0 algebry, realizuemyh
operatorami Gilp1\-ber\-ta-Shmidta pri tochnom predstavlenii v modulyah Verma
nad algebroe0 Li $\sltwo$ dlya vseh dostatochno malyh znachenie0 $\hbar$
(inymi slovami, esli $\pi_{h}$ -- predstavlenie $\frak A$ v $V_h$ i $\HS_h$ --
prostranstvo operatorov Gilp1berta-Shmidta v $V_h$, to $\HS$ ponimaet\-sya
v teoreme v smysle
$$\lim_{\hbar\to0}\bigcap_{h_0<h<h_0+\hbar}\pi_{h}^{-1}(\HS_{h}\cap
\pi_{h}(\frak A)),$$
gde $h_0=\tfrac12$ ili $h_0=1$).

\demo\nofrills{}\cyr Rezulp1tat teoremy 3B sleduet iz yavnyh vychislenie0.
Privedem vy\-chis\-le\-nie kommutatora $[L_2,L_{-2}]$ i ego otsenki, kotoryh
dostatochno dlya op\-re\-de\-le\-niya tsentralp1nogo zaryada algebry Virasoro.
Modulp1 Verma $V_h$ realizovan v prostranstve mnogochlenov $\Bbb C[z]$,
$\xi=z\partial_z$, a generatory $q_R$--konformnyh simmetrie0 zadayut\-sya
vyrazheniyami, vypisannymi ranee. Tog\-da pri $h\to\frac12$
$$\aligned
[L_2,L_{-2}]=&\tfrac{(\xi+3h)^2(\xi+1)(\xi+2)}{(\xi+2h)(\xi+2h+1)}-
\tfrac{(\xi+3h-2)^2\xi(\xi-1)}{(\xi+2h-1)(\xi+2h-2)}=\\
&\tfrac{(\xi+\tfrac32+3\hbar)^2(\xi+1)(\xi+2)}{(\xi+1+2\hbar)(\xi+2+2\hbar)}-
\tfrac{(\xi-\tfrac12+3\hbar)^2\xi(\xi-1)}{(\xi-1+2\hbar)(\xi+2\hbar)}.
\endaligned$$
Imeyet mesto sleduyushchaya otsenka ``$\md O(\hbar^2)$'' ($\hbar$ v
znamenatelyah ostavlen, t.k. pri nulevom $\hbar$ ne opredeleno
predstavlenie rasshirennoe0 algebry Vee0lya v sootvet\-stvuyushchem module
Verma):
$$\aligned
[L_2,L_{-2}]\equiv &(\xi+\tfrac32+3\hbar)^2(1-\tfrac{2\hbar}{\xi+1+\hbar})
(1-\tfrac{2\hbar}{\xi+2+\hbar})-\\
&(\xi-\tfrac12+3\hbar)^2(1-\tfrac{2\hbar}{\xi-1+\hbar})
(1-\tfrac{2\hbar}{\xi+\hbar})\mod{O(\hbar^2)}.
\endaligned$$
Vydelyaya chleny poryadka $\hbar$ (poskolp1ku starshie0 chlen poryadka
edinitsy uzhe izvesten), poluchim:
$$2\hbar(\xi-\tfrac12)^2(\tfrac1{\xi-1+\hbar}+\tfrac1{\xi+\hbar})-
2\hbar(\xi+\tfrac32)^2(\tfrac1{\xi+1+\hbar}+\tfrac1{\xi+2+\hbar})+12\hbar.$$
Provedem teperp1 otsenku ``$\md\HS$'' i poluchim:
$$\aligned
&2\hbar\xi^2(\tfrac1{\xi-1+\hbar}+\tfrac1{\xi+\hbar}-\tfrac1{\xi+1+\hbar}-
\tfrac1{\xi+2+\hbar})-12\hbar-\\
&2\hbar\xi(\tfrac1{\xi-1-\hbar}+\tfrac1{\xi+\hbar})-
6\hbar(\tfrac1{\xi+1+\hbar}+\tfrac1{\xi+2+\hbar})\sim\\
&2\hbar\xi^2(\tfrac2{(\xi-1+\hbar)(\xi+1+\hbar)}+
\tfrac2{(\xi+\hbar)(\xi+2+\hbar)})+12\hbar-16\hbar\sim\\
&8\hbar+12\hbar-16\hbar=4\hbar.
\endaligned$$
Takim obrazom, ``asimptoticheskie0'' tsentralp1nye0 zaryad raven $8\hbar$.

Analogichno, pri $h\to1$
$$\aligned
[L_2,L_{-2}]=&\tfrac{(\xi+3h)^2(\xi+1)(\xi+2)}{(\xi+2h)(\xi+2h+1)}-
\tfrac{(\xi+3h-2)^2\xi(\xi-1)}{(\xi+2h-1)(\xi+2h-2)}=\\
&\tfrac{(\xi+3+3\hbar)^2(\xi+1)(\xi+2)}{(\xi+2+2\hbar)(\xi+3+2\hbar)}-
\tfrac{(\xi+1+3\hbar)^2\xi(\xi-1)}{(\xi+2\hbar)(\xi+1+2\hbar)}.
\endaligned$$
Imeet mesto sleduyushchaya otsenka ``$\md O(\hbar^2)$'' ($\hbar$ v
znamenatelyah ostavlen, t.k. pri nulevom $\hbar$ ne opredeleno
predstavlenie rasshirennoe0 algebry Vee0lya v sootvet\-stvuyushchem
module Verma):
$$\aligned
[L_2,L_{-2}]\equiv &(\xi+1)(\xi+3+3\hbar)(1-\tfrac{2\hbar}{\xi+2+\hbar})
(1-\tfrac{2\hbar}{\xi+3+\hbar})(1+\tfrac{3\hbar}{\xi+3+\hbar})-\\
&(\xi-1)(\xi+1+3\hbar)(1-\tfrac{2\hbar}{\xi+\hbar})
(1-\tfrac{2\hbar}{\xi+1+\hbar})(1+\tfrac{3\hbar}{\xi+1+\hbar})\mod{O(\hbar^2)}.
\endaligned$$
Vydelyaya chleny poryadka $\hbar$ (poskolp1ku starshie0 chlen poryadka
edinitsy uzhe izvesten), poluchim:
$$\hbar(\xi+1)(\xi+3)(\tfrac1{\xi+3+\hbar}-\tfrac2{\xi+2+\hbar})-
\hbar(\xi+1)(\xi-1)(\tfrac1{\xi+1+\hbar}-\tfrac2{\xi+\hbar})+6\hbar.$$
Provedem teperp1 otsenku ``$\md\HS$'' i poluchim:
$$8\hbar+2(\xi+1)(\tfrac{\xi-1}{\xi+\hbar}-\tfrac{\xi+3}{\xi+2+\hbar})\sim
8\hbar-4\hbar=4\hbar.$$
Takim obrazom, ``asimptoticheskie0'' tsentralp1nye0 zaryad opyatp1 raven
$8\hbar$.
\enddemo

\remark{\cyi Zamechanie 9}\cyr Otsenka ``$\md O(\hbar^2)$'' mozhet bytp1
uluchshena do ``$\md O(\hbar^n)$'' dlya lyubogo {\cyi konechnogo\/} $n$ (no
ne do ``$\md O(\hbar^\infty)$'' iz-za posleduyushchee0 otsenki ``$\md\HS$'').
\endremark

\remark{\cyi Zamechanie 10}\cyr Kak sleduet iz rezulp1tatov vtorogo
paragrafa pe\-res\-ta\-nov\-ka poryadka otsenok privodit k nulevomu
tsentralp1nomu zaryadu dlya algebry Virasoro.
\endremark

\cyr Po-vidimomu, neperestanovochnostp1 poryadka otsenok ``$\md O(\hbar^n)$''
i ``$\md\mathbreak\HS$'' v sluchae $q_R$--konformnyh simmetrie0 yavlyaet\-sya
otrazheniem bolee ob\-shche\-go i fundamentalp1nogo fakta razlichie0 mezhdu
asimptoticheskoe0 te\-o\-ri\-ee0 psevdodifferentsialp1nyh operatorov i
psevdodifferentsialp1nym ischisleniem na asimptoticheskih mnogoobraziyah
v smysle [27].

Otmetim, chto rezulp1tat teoremy 3B soderzhit strannoe sovpadenie
``asimptoticheskih'' tsentralp1nyh zaryadov dlya algebry Virasoro pri
$h\to\frac12$ i $h\to1$ (a imenno, $c=8\hbar$, $\hbar\to0$). Vozmozhno,
e1to sovpadenie obp2yasnyaet\-sya tem, chto pri vseh znacheniyah
e1kstremalp1nogo vesa $h$ algebra Virasoro prisut\-stvuet na nekotorom
eshche bolee skrytom i poka ne vyyavlennom urovne universalp1no (t.e. ee
harakteristiki ne zavisyat ot
$h$).

\head\cyb Zaklyuchenie
\endhead

V rabote issledovany razlichnye priblizhennye predstavleniya
bes\-ko\-nech\-no\-mer\-nyh $\Bbb Z$--graduirovannyh algebr Li: algebry
Vitta loranovskih polinomialp1nyh vektornyh polee0 na okruzhnosti i ee
odnomernogo ne\-tri\-vi\-alp1\-no\-go tsentralp1nogo rasshireniya, algebry
Virasoro, bes\-ko\-nech\-no\-mer\-ny\-mi skrytymi simmetriyami v modulyah
Verma nad algebroe0 Li $\sltwo$. Rassmotreny kak asimptoticheskie
predstavleniya ``$\md O(\hbar^n)$'' i predstavleniya s tochnostp1yu do
operatorov iz nekotorogo klassa $\frak S$ (kom\-pak\-t\-nyh operatorov,
operatorov Gilp1berta-Shmidta ili operatorov ko\-nech\-no\-go ranga), tak i
sluchai, sovmeshchayushchie oba tipa priblizhenie0 (i v dannom sluchae
vyyavlen e1ffekt neperestanovochnosti poryadka ih vy\-pol\-ne\-niya, chto,
po-vidimomu, svidetelp1stvuet o bolee obshchem i fun\-da\-men\-talp1\-nom
fakte razlichie0 mezhdu asimptoticheskoe0 teoriee0
psev\-do\-dif\-fe\-ren\-tsi\-alp1\-nyh operatorov i psevdodifferentsialp1nym
ischisleniem na asimptoticheskih mnogoobraziyah v smysle [27]). Nekotorye
prilozheniya obsuzhdaemyh voprosov k prikladnym problemam informatsionnyh
teh\-no\-lo\-gie0 (organizatsiya peredachi informatsii v integrirovannyh
vi\-deo\-kog\-ni\-tiv\-nyh interaktivnyh sistemah dlya uskorennyh
neverbalp1nyh kom\-pp1yu\-ter\-nyh i telekommunikatsie0) izuchalisp1 v
rabote [28].

\rm\ \newline
\ \newline
\centerline{\cybe Spisok literatury}
\eightpoint
\roster
\item" [1]" {\cyre Te1e0lor M., Psevdodifferentsialp1nye operatory. M., Mir,
1985.}
\item" [2]" {\cyre Trev F., Vvedenie v psevdodifferentsialp1nye operatory
i integralp1nye ope\-ra\-to\-ry Furp1e. M., Mir, 1984.}
\item" [3]" {\cyre Bogolyubov N.N., Mitropolp1skie0 Yu.A., Asimptoticheskie
metody v teorii nelinee0nyh kolebanie0. M., Fizmatlit, 1958.}
\item" [4]" {\cyre Maslov V.P., Teoriya vozmushchenie0 i asimptoticheskie
metody. M., Izd-vo MGU, 1965.}
\item" [5]" {\cyre Maslov V.P., Fedoryuk M.V., Kvaziklassicheskoe priblizhenie
dlya uravnenie0 kvantovoe0 mehaniki. M., Nauka, 1976.}
\item" [6]" {\cyre Karasev M.V., Maslov V.P., Nelinee0nye skobki Puassona.
Geometriya i kvan\-to\-va\-nie. M., Nauka, 1991.}
\item" [7]" {\cyre Kirillov A.A., E1lementy teorii predstavlenie0. M., Nauka,
1978.}
\item" [8]" {\cyre Barut A., Ronchka R., Teoriya predstavlenie0 grupp i ee
prilozheniya. M., Mir, 1980.}
\item" [9]" Juriev D., Topics in hidden symmetries. V. E-print:
funct-an/9611003 (1996).
\item"[10]" Juriev D., On the infinite dimensional hidden symmetries. III.
$q_R$--conformal symmetries at $q_R\!\to\!\infty$ and Berezin-Karasev-Maslov
asymptotic quantization of $C^\infty(S^1)$. E-print: funct-an/9702002 (1997).
\item"[11]" {\cyre Zhelobenko D.P., Predstavleniya reduktivnyh algebr Li.
M., Nauka, 1993.}
\item"[12]" {\cyre Diksmp1e Zh., Universalp1nye obertyvayushchie algebry.
M., Mir, 1976.}
\item"[13]" {\cyre Yurp1ev D., Kompleksnaya proektivnaya geometriya i
kvantovaya proektivnaya te\-o\-riya polya /\!/ TMF. 1994. T.101. vyp.3.
S.331-348}.
\item"[14]" {\cyre Berezin F.A., Kvantovanie v kompleksnyh simmetricheskih
prostranstvah /\!/ Izvestiya AN SSSR. Ser.matem. 1975. T.39, vyp.2. S.363-402}.
\item"[15]" Klimek S., Lesniewski A., Quantum Riemann surfaces. I /\!/
Commun.Math.Phys. 1992. V.146. P.103-122.
\item"[16]" {\cyre Bidenharn L., Lauk Dzh., Uglovoe0 moment v kvantovoe0
fizike. M., Mir, 1984}.
\item"[17]" {\cyre Fuks D.B., Kogomologii beskonechnomernyh algebr Li. M.,
Nauka, 1984.}
\item"[18]" Juriev D., Noncommutative geometry, chiral anomaly in the quantum
projective ($\sltwo$--invariant) field theory and $\jltwo$--invariance /\!/
J.Math.Phys. 1992. V.33. P.2819-2822, (E): 1993. V.34. P.1615.
\item"[19]" Juriev D., Remarks on nonassociative structures in quantum
projective field theory: the central extension $\jltwo$ of the double
$\sltwo+\sltwo$ of the simple Lie algebra $\sltwo$ and related topics.
Acta Appl.Math. 1998. V.50. P.191-196.
\item"[20]" {\cyre Yurp1ev D.V., Kvantovaya proektivnaya teoriya polya:
kvantovo-polevye analogi uravnenie0 E1e0lera-Arnolp1da v proektivnyh
$G$--gipermulp1tipletah /\!/ TMF. 1994. T.98, vyp.2. S.220-240}.
\item"[21]" {\cyre Malp1tsev A.I., O predstavlenii neassotsiativnyh kolets /\!/
UMN. 1952. T.7, vyp.1. S.181-185 [pereizdano v ``Izbrannyh trudah. 1.
Klassicheskaya al\-geb\-ra'', M., Nauka, 1976, S.328-331]}.
\item"[22]" {\cyre Kobayasi Sh., Nomidzu K., Osnovy differentsialp1noe0
geometrii. M., Nauka, 1981}.
\item"[23]" Juriev D., Infinite dimensional geometry and quantum field theory
of strings. II. Infinite-dimensional noncommutative geo\-met\-ry of a
self-interacting string field /\!/ Russian J.Math. Phys. 1996. V.4, no.3.
P.287-314.
\item"[24]" Juriev D., On the infinite dimensional hidden symmetries. I.
Infinite dimensional geo\-met\-ry of $q_R$--conformal symmetries.
E-print: funct-an/9612004 (1996).
\item"[25]" Juriev D., On the infinite dimensional hidden symmetries. II.
$q_R$--conformal modular functor. E-print: funct-an/9701009 (1997).
\item"[26]" {\cyre Kirillov A.A., Vvedenie v teoriyu predstavlenie0 i
nekommutativnye0 gar\-mo\-ni\-ches\-kie0 analiz / Sovr.probl.mat.
Fund.napravleniya. T.22. M., VINITI, 1987, S.5-162}.
\item"[27]" Dubnov V.L., Maslov V.P., Naza{\v\i}kinski{\v\i} V.E.,
The complex Lagrangian germ and the canonical operator /\!/ Russian
J.Math.Phys. 1995. V.3. P.141-180.
\item"[28]" {\cyre Yurp1ev D.V., Droe1my: e1ksperimentalp1naya matematika,
informatika i bes\-ko\-nech\-no\-mer\-naya geometriya:} Report RCMPI-96/05$^+$
(1996).
\endroster
\enddocument